\documentclass{amsart}
\usepackage{fullpage}
\usepackage[all]{xy}
\usepackage{graphics}
\input{pramod.tex}

\CompileMatrices

\newenvironment{genthm}[1]{\begingroup\theoremstyle{plain}%
  \newtheorem*{mygenthm}{#1}\begin{mygenthm}}{\end{mygenthm}\endgroup}

\newcommand{\PP}{\mathcal{P}}
\newcommand{\Pe}{\PP_\epsilon}
\newcommand{\tP}{\tilde\PP}
\newcommand{\tPsp}{\tilde\PP^{\text{\rm sp}}}
\newcommand{\tPe}{\tilde\PP_\epsilon}
\newcommand{\tPo}{\tilde\PP^\circ}

\newcommand{\tPspe}{\tilde\PP^{\text{\rm sp}}_\epsilon}
\renewcommand{\N}{\mathcal{N}}
\newcommand{\No}{\mathcal{N}_{\text{\rm o}}}
\newcommand{\Nosp}{\No^{\text{\rm sp}}}
\newcommand{\Noc}{\mathcal{N}_{\text{\rm o,c}}}
\newcommand{\Nocc}{\mathcal{N}_{\text{\rm o,\=c}}}
\newcommand{\Noccsp}{\mathcal{N}_{\text{\rm o,\=c}}^{\text{\rm sp}}}
\newcommand{\Nor}{\mathcal{N}_{\text{\rm o,r}}}
\newcommand{\Noro}{\mathcal{N}^0_{\text{\rm o,r}}}
\newcommand{\Norro}{\mathcal{N}^0_{\text{\rm o,\=r}}}
\newcommand{\vn}{\varnothing}

\newcommand{\Ab}{{\bar A}}
\newcommand{\pr}[1]{\mathit{pr}_{#1}}
\newcommand{\dLS}{d_{\mathrm{LS}}}
\newcommand{\dBV}{d_{\mathrm{BV}}}
\newcommand{\dS}{d_{\mathrm{S}}}
\newcommand{\dA}{\bar d}
\newcommand{\db}{\bar d}
\newcommand{\cmp}[2]{{}^{\langle{#2}\rangle\!}{#1}}
\newcommand{\pmp}[2]{\cmp{#1}{#2}}

\title[An order-reversing duality map]{An order-reversing duality map
for conjugacy classes in Lusztig's canonical quotient}

\author{Pramod N. Achar}
\date{16 September 2002}
\address{Department of Mathematics\\
  University of Chicago\\
  Chicago, IL \ 60637}
\email{pramod@math.uchicago.edu}

\begin{document}

\begin{abstract}
We define a partial order on the set $\Nocc$ of pairs $(\orb,C)$,
where $\orb$ is a nilpotent orbit and $C$ is a conjugacy class in
$\Ab(\orb)$, Lusztig's canonical quotient of $A(\orb)$.  We then show
that there is a unique order-reversing duality map $\Nocc \to
\Langdual\Nocc$ that has certain properties analogous to those
of the original Lusztig-Spaltenstein duality map.  This generalizes
work of E.~Sommers.
\end{abstract}

\maketitle

\section{Introduction}
\label{sect:intro}

Let $G$ be a connected simple complex algebraic group, and let
$\Lie{g}$ be its Lie algebra.  Let $\N$ be the nilpotent cone in
$\Lie{g}$; let $\No$ be the set of $G$-orbits in $\N$.  The notion of
a duality map for nilpotent orbits has its roots in the introduction
of ``special representations'' of a Weyl group by Lusztig
\cite{lusztig:special}.  He gave a purely algebraic treatment aimed at
studying primitive ideals in enveloping algebras, but in passing, he
conjectured that (and it was quickly verified that) all special
representations should be assigned to nilpotent orbits with the
trivial local system via Springer's correspondence, so that special
representations would sit in bijection with a remarkable set of
``special'' nilpotent orbits, denoted $\Nosp$.  Subsequently, Lusztig
and Spaltenstein observed that the set of special nilpotent orbits
admits a natural order-reversing bijection (with respect to the usual
closure order on nilpotent orbits) that usually corresponds, in the
language of special representations, to tensoring with the sign
representation.  (There are a couple of curious exceptions to this in
types $E_7$ and $E_8$.)  Indeed, this bijection could be extended to
an order-reversing map $\dLS: \No \to \No$ whose image consists
precisely of the special orbits, and which is an involution when
restricted to its image.  In \cite{spaltenstein:classes}, Spaltenstein
gives an axiomatic treatment of the map $\dLS$, showing that it is the
unique map satisfying certain order conditions and a certain
compatibility with induction.

Now, since a group $G$ and its Langlands dual $\Langdual G$ have
isomorphic Weyl groups, there is a natural bijection between their
respective sets of special nilpotent orbits $\Nosp$ and
$\Langdual\Nosp$.  Spaltenstein observed that this bijection is
order-preserving, so by composing $\dLS$ with it, one obtains a map
$\No \to \Langdual\No$ or $\Langdual\No \to \No$.  Barbasch and Vogan
later gave an elegant and intrinsic construction of this incarnation
of the map, which we shall denote by $\dBV$, in terms of associated
varieties of certain Harish-Chandra modules.

Sommers \cite{sommers:duality} has shown how to enlarge the domain of
$\dBV$ so that the extended map surjects onto $\Langdual\No$.  This
latter set does not, in general, sit in bijection with $\No$, so there
is no analogue of Sommers' map for $\dLS$.  Let $A(\orb)$ be the
component group of the centralizer in $G$ of some element of $\orb$,
and let $\Noc$ be the set of pairs $(\orb,C)$, where $\orb \in \No$
and $C$ is a conjugacy class in $A(\orb)$.  (We do not need to be
careful about which element of $\orb$ we pick to define $A(\orb)$,
since any two yield component groups that are canonically isomorphic
up to inner automorphism.)  Sommers' map $\dS: \Noc \to \Langdual\No$
agrees with $\dBV$ when composed with the inclusion $\No
\hookrightarrow \Noc$ defined by $\orb \mapsto (\orb,1)$.

Finally, let $\Ab(\orb)$ be Lusztig's canonical quotient of $A(\orb)$.
This was originally introduced by Lusztig \cite{lusztig:characters}
for special orbits, but Sommers \cite{sommers:duality}, in the course
of giving a new characterization of the canonical quotient, observes
that the definition makes sense for all orbits.  We let $\Nocc$ be the
set of pairs $(\orb,C)$, where this time $C$ is a conjugacy class in
$\Ab(\orb)$.  Sommers' description of $\Ab(\orb)$ leads to a
proof of the following statement (\cite{sommers:duality},
Proposition~15): if $C$ and $C'$ are two conjugacy classes in
$A(\orb)$ that descend to the same conjugacy class in $\Ab(\orb)$,
then $\dS(\orb,C) = \dS(\orb,C')$.  In other words, $\dS$ factors
through the natural projection $\Noc \twoheadrightarrow \Nocc$.  In
this article, we often regard $\dS$ as a map $\Nocc \to \Langdual\No$.

One task we accomplish in this paper is the introduction of a partial
order on the set $\Nocc$, as follows.  We say that $(\orb,C) \le
(\orb',C')$ if
\begin{equation}\label{eqn:po}
\orb \le \orb'
\qquad\text{and}\qquad
\dS(\orb,C) \ge \dS(\orb',C').
\end{equation}
\emph{A priori}, this partial order might not be well-defined: we
might have had $\dS(\orb,C) = \dS(\orb,C')$ even when $C \ne C'$.  In
the course of this paper, we rectify this by proving a converse to
Proposition~15 of \cite{sommers:duality}.

\begin{thm-}\label{thm:po}
Let $C, C' \subset A(\orb)$ be two conjugacy classes associated to the
same orbit.  Then $\dS(\orb,C) = \dS(\orb,C')$ if and only if $C$ and
$C'$ have the same image in $\Ab(\orb)$.  As a consequence, the
partial order (\ref{eqn:po}) on $\Nocc$ is well-defined.
\end{thm-}

The principal aim of this paper is to show that $\Nocc$ admits a
unique duality map $\db$ that is compatible with the aforementioned
maps in the appropriate senses.  In particular, such a duality map
ought to satisfy the partial-order properties of $\dLS$ and $\dBV$:
\begin{enumerate}
\item \label{ax:ord-rev} If $(\orb,C) \le (\orb',C')$, then
$\db(\orb,C) \ge \db(\orb',C')$. 
\item \label{ax:ord-pres} $\db^2(\orb,C) \ge (\orb,C)$.
\end{enumerate}
It also ought to coincide with $\dBV$ and $\dS$ when its domain or
codomain is restricted.  Indeed, we need only make an explicit
requirement with respect to $\dS$; that automatically implies the
desired compatibility with $\dBV$ as well.  We write $\pr1: \Nocc \to
\No$ for the obvious projection.
\begin{enumerate}
\setcounter{enumi}{2}
\item \label{ax:dS-resp} $\pr1 \circ \db(\orb,C) = \dS(\orb,C)$.
\end{enumerate}
Finally, we need one additional condition to guarantee the uniqueness
of the map.
\begin{enumerate}
\setcounter{enumi}{3}
\item \label{ax:im-size} Among maps respecting the first three axioms,
$\db$ has an image set of maximal size.
\end{enumerate}

\noindent
To be precise, we ought to be seeking a pair of maps $\db :
\Nocc \to \Langdual\Nocc$, $\db: \Langdual\Nocc \to \Nocc$,
both of which satisfy the above axioms.  Indeed,
axiom~(\ref{ax:ord-pres}) only makes sense if we have two such maps
together.  Nevertheless, to avoid making the language too cumbersome,
we will speak throughout the paper of ``a'' duality map $\Nocc \to
\Langdual\Nocc$, and always assume it to be implicitly accompanied by
a partner map $\Langdual\Nocc \to \Nocc$.

The main result of the paper is the following.

\begin{thm-}\label{thm:d}
There is a unique map $\db: \Nocc \to \Langdual\Nocc$ satisfying the
axioms (\ref{ax:ord-rev})--(\ref{ax:im-size}).
\end{thm-}

\noindent
Let us call this map the \emph{extended duality map}.  

In type $A$, of course, all the $\Ab(\orb)$-groups are trivial, so
this theorem does not say anything new: the extended duality map is
just the same as $\dBV$.  In all other types, the theorem will be
proved by giving an explicit construction of the map.  For the
classical types, this entails a combinatorial algorithm in terms of
partitions, whereas in the exceptional groups, we define $\db$ simply
by tabulating all its values.

We begin our discussion in Section~\ref{sect:formal} by collecting
some properties that must be satisfied by any putative extended
duality map.  These lead up to a criterion for showing that a
candidate map satisfies axiom~(\ref{ax:im-size}), and that it is the
unique such map.  In Section~\ref{sect:partitions}, we define the
combinatorial objects that will be used to work with $\No$ and $\Nocc$
in the classical groups, and we recall various useful facts about
them.  Serious work on the classical-groups case begins in
Section~\ref{sect:constr-class}, where we give the definitions of
$\dA$ and develop some basic techniques for studying it.
Section~\ref{sect:proof-class} contains the proof of
Theorems~\ref{thm:po} and~\ref{thm:d} for the classical groups.  In
Section~\ref{sect:explicit}, we consider the exceptional groups, for
which the main theorems are proved simply by drawing out the
partial-order diagrams and verifying the existence of the extended
duality map by inspection.  Additionally, these partial-order diagrams
are accompanied by those for a few classical groups of low rank,
simply for the sake of having some examples to look at.  Finally, in
Section~\ref{sect:comments}, we explore some possible applications and
consequences of the present work.

I would like to thank A.-M.~Aubert, R.~Bezrukavnikov, V.~Ginzburg,
R.~Kottwitz, V.~Ostrik, and D.~Vogan for helpful conversations.  I
would like to specifically thank one of my referees for proposing
axiom~(\ref{ax:im-size}): an earlier draft of this work employed a
different statement, which did not imply uniqueness.  Finally, I would
especially like to thank E.~Sommers.  His paper \cite{sommers:duality}
provides the bulk of the motivation for this one; this paper would not
have been possible in the absence of the numerous discussions I have
had with him on these topics.

\section{Formal properties of duality}
\label{sect:formal}

Throughout this section, we assume that Theorem~\ref{thm:po} is true,
so the partial order on $\Nocc$ is defined, and the axioms for an
extended duality map make sense.  We begin a few easy properties of
this partial order.

\begin{prop}\label{prop:inherit-po}
We have $\orb \le \orb'$ in $\No$ if and only if $(\orb,1) \le
(\orb',1)$ in $\Nocc$.  Thus, via the imbedding $\orb \mapsto
(\orb,1)$, the partially ordered set $\No$ can be regarded as a subset
of $\Nocc$ with the inherited partial order.
\end{prop}
\begin{proof}
From the definition of the partial order, we know that $(\orb,1) \le
(\orb',1)$ implies that $\orb \le \orb'$.  For the converse, we need
to prove that if $\orb \le \orb'$, then $\dS(\orb,1) \ge
\dS(\orb',1)$.  But we know that $\dS(\orb,1) = \dBV(\orb)$ and
$\dS(\orb',1) = \dBV(\orb')$, and we further know that $\orb \le
\orb'$ implies $\dBV(\orb) \ge \dBV(\orb')$.
\end{proof}

\begin{prop}
Regarding the Sommers duality map $\dS$ as being a map $\Nocc \to
\Langdual\No$, we have that $(\orb,C) \le (\orb',C')$ implies
$\dS(\orb,C) \ge \dS(\orb',C')$.  That is, $\dS$ is an order-reversing
map.
\end{prop}
\begin{proof}
This is an obvious consequence of the definition of the partial order
on $\Nocc$.
\end{proof}

\begin{prop}\label{prop:triv-min}
For a fixed orbit $\orb$ and any conjugacy class $C \subset
\Ab(\orb)$, we have $(\orb,1) \le (\orb,C)$.
\end{prop}
\begin{proof}
(This fact is hinted at in \cite{sommers:duality}, where it is proved
that $(\orb,1)$ has minimal $\tilde{b}$-value among all the
$(\orb,C)$.)  All we have to check is that $\dS(\orb,1) \ge
\dS(\orb,C)$.  In the exceptional groups, we can verify this simply by
scanning Sommers' tables of computed values from
\cite{sommers:duality}.  In the classical groups, it is an easy
computation from Sommers' formulas for $\dS$, which we recall at the
end of Section~\ref{sect:partitions}.  We defer carrying out the
computation until then.
\end{proof}

We now turn our attention to duality maps.  Some formal properties can
be deduced from just the first three axioms.  Let us define a
\emph{weak extended duality map} to be any map $\db: \Nocc \to
\Langdual\Nocc$ satisfying the first three axioms, but not
necessarily the fourth.  Let us say that a pair $(\orb,C)$ is
\emph{special for $\db$}, or simply \emph{special} if no ambiguity is
likely, if it is in the image of $\db$.

\begin{prop}\label{prop:d3-d}
We have $\db^3 = \db$, so that when we restrict to the special set,
the map $\db$ is an order-reversing bijection between special pairs in
$\Nocc$ and those in $\Langdual\Nocc$, and $\db^2$ is the identity
map.  In general, $\db^2(\orb,C)$ is the unique smallest special
pair that is greater than or equal to $(\orb,C)$.
\end{prop}
\begin{proof}
Axiom~(\ref{ax:ord-pres}) says that $\db^2(\orb,C) \ge (\orb,C)$.
Applying $\db$ to both sides of this, we obtain $\db^3(\orb,C) \le
\db(\orb,C)$, by axiom~(\ref{ax:ord-rev}).  But on the other hand,
axiom~(\ref{ax:ord-pres}) also tells us that $\db^2(\db(\orb,C)) \ge
\db(\orb,C)$.  We conclude that $\db^3(\orb,C) = \db(\orb,C)$.

For the second part of the proposition, we know that $\db^2(\orb,C)$
is special and greater than or equal to $(\orb,C)$.  Now, let
$(\orb',C') \ge (\orb,C)$ be any special pair.  We have $\db(\orb',C')
\le \db(\orb,C)$, whence $\db^2(\orb',C') \ge \db^2(\orb,C)$.  
But since $(\orb',C')$ is special, we have $\db^2(\orb',C') =
(\orb',C')$, so we can deduce that $(\orb',C') \ge \db^2(\orb,C)$:
thus $\db^2(\orb,C)$ is the smallest special pair that is greater
than or equal to $(\orb,C)$.
\end{proof}

\begin{prop}\label{prop:spec-uniq}
If $\db_1, \db_2: \Nocc \to \Langdual\Nocc$ are two weak extended
duality maps giving rise to the same special set, then $\db_1 =
\db_2$.
\end{prop}
\begin{proof}
We first show that $\db_1$ and $\db_2$ agree on special pairs.  Suppose
$(\orb,C)$ is special, and that $\orb' = \dS(\orb,C)$.  We must have
$\db_1(\orb,C) = (\orb',C_1)$, $\db_2(\orb,C) = (\orb',C_2)$ for some
$C_1$ and $C_2$.  Moreover, the preceding proposition tells us that
$\db_1(\orb',C_1) = (\orb,C) = \db_2(\orb',C_2)$.  Applying
axiom~(\ref{ax:dS-resp}) again, we have $\dS(\orb',C_1) =
\dS(\orb',C_2) = \orb$.  Finally, Theorem~\ref{thm:po} says that we
must have $C_1 = C_2$.

Second, if $(\orb,C)$ is nonspecial, Proposition~\ref{prop:d3-d} tells
us that there exists a unique smallest special pair $(\orb_0,C_0)$
that is larger than $(\orb,C)$, and that $\db_1(\orb,C) =
\db_1(\orb_0,C_0) = \db_2(\orb_0,C_0) = \db_2(\orb,C)$.
\end{proof}

Within the proof of this last proposition lurks an important
observation: all weak extended duality maps $\db$ for which a given
pair $(\orb,C)$ is special take the same value on it.  This is
because, by Theorem~\ref{thm:po}, there is at most one class $C'$ such
that $\dS(\orb',C') = \orb$, where $\orb' = \dS(\orb,C)$.  If there
does not exist such a $C'$, then $(\orb,C)$ cannot be special for any
weak extended duality map.  Inspired by this, we define the set
\[
\Noccsp = \{ (\orb,C) \in \Nocc \mid \text{there exists a $C'$ such
that $\dS(\orb',C') = \orb$, where $\orb' = \dS(\orb,C)$} \},
\]
and note that the special set of any weak extended duality map must be
contained within $\Noccsp$.  If there exists one whose special set is
the entirety of $\Noccsp$, then it would automatically satisfy the
fourth axiom as well.  It would also be the unique possible extended
duality map, by Proposition~\ref{prop:spec-uniq}.  We have established
the following.

\begin{prop}\label{prop:exist-uniq}
If $\db$ is a weak extended duality map whose special set is
$\Noccsp$, then $\db$ is in fact the unique extended duality map. \qed
\end{prop}

\noindent
The proof of Theorem~\ref{thm:d} in Sections~\ref{sect:proof-class}
and~\ref{sect:explicit} is carried out by explicitly constructing a
weak extended duality map $\db$ which happens to have all of $\Noccsp$
as its special set, and then applying
the preceding proposition.  Although that construction is a
laborious undertaking which occupies most of this paper, there is
\emph{a posteriori} a concise, uniform description of the extended
duality map.  It is the map whose image is $\Langdual\Noccsp$, and
whose values are computed according to the discussion in the proof of
Proposition~\ref{prop:spec-uniq}.

We conclude this section with a few additional observations about the
extended duality map.  Below, $\db$ will denote only the extended
duality map, and ``special'' will refer to all elements of $\Noccsp$.

\begin{prop}\label{prop:orbit-special}
Any pair of the form $(\orb,1)$ is special.
\end{prop}
\begin{proof}
Let $\orb' = \dS(\orb,1) = \dBV(\orb)$.  To show that $(\orb,1) \in
\Noccsp$, we must merely demonstrate the existence of a class $C'$
such that $\dS(\orb',C') = \orb$.  This was done by Sommers in
\cite{sommers:duality} with his construction of a ``canonical
inverse'': this is a certain right inverse to $\dS$ that was used to
show that $\dS$ is surjective.  The details of the construction are
such that the preimage produced for $\orb$ is always a conjugacy class
associated to $\dBV(\orb)$.
\end{proof}

The following two statements are easily deduced from the above
uniform description of $\db$.

\begin{prop}\label{prop:dS-inv}
Sommers' canonical inverse is given by $\orb \mapsto
\db(\orb,1)$. \qed
\end{prop}

\begin{prop}\label{prop:special-orbits}
An orbit $\orb$ is special if and only if $\db(\orb,1) =
(\dBV(\orb),1)$. \qed
\end{prop}

\noindent
Even if $\orb$ is a special orbit, we cannot say anything in general
about whether $(\orb,C)$ is a special pair for nontrivial $C$.  The
computed examples in Section~\ref{sect:explicit} include instances of
both special and nonspecial pairs of this form.

\section{Orbits, partitions, and component groups}
\label{sect:partitions}

We spend this section collecting facts and formulas for working with
partitions as a way of understanding nilpotent orbits in the classical
groups.  It is suggested that the reader skip this section, referring
back to it only when necessary to find a particular definition or
formula.

\subsection{Partitions}

Let $\PP(n)$ be the set of partitions of $n$.  For a partition
$\lambda$, let $|\lambda|$ denote the sum of the parts of $\lambda$.
We typically write $\lambda = [ \lambda_1 \ge \lambda_2 \ge \cdots \ge
\lambda_k ]$, and we assume $\lambda_k \ne 0$ unless stated otherwise.
Sometimes, however, we shall write partitions as follows, using
exponents to indicate multiplicities: $[a_1^{p_1}, \ldots,
a_k^{p_k}]$, with $a_1 > \cdots > a_k$.  Let $r_\lambda(a)$, or simply
$r(a)$, denote the multiplicity of $a$ as a part in $\lambda$.  We
define the \emph{height} of a part in a partition to be the number of
parts greater than or equal to the given one:
$
\height_\lambda(a) = \height(a) = \sum_{b \ge a} r_\lambda(b)
$.
Note that this formula makes sense even if $a$ is not a part of
$\lambda$; {\it i.e.}, if $r_\lambda(a) = 0$.  We shall employ the
notion of height in such circumstances from time to time; we may refer
to it as ``generalized height'' to draw attention to the fact that
$r_\lambda(a) = 0$.  Finally, we write $\#\lambda$ to denote the total
number of parts of $\lambda$.

For odd $n$, we write $\PP_B(n)$ for the set of partitions in which
even parts occur with even multiplicity.  For even $n$, we write
$\PP_C(n)$ for the set of partitions in which odd parts occur with
even multiplicity, and $\PP_D(n)$ for the set of partitions in which
even parts occur with even multiplicity.  Here, the subscript letters
correspond to the type of classical Lie group whose nilpotent orbits
are indexed by the given set of partitions, with one caveat: very even
partitions (those consisting only of even parts with even
multiplicity) in type $D$ correspond to two nilpotent orbits.  We
ignore this fact throughout the paper, because such orbits have
trivial $A(\orb)$-groups, so the duality map we construct here will
not have anything new to say about them.  We will sometimes write
$\PP_1(n)$ for $\PP_C(n)$, and $\PP_0(n)$ for either $\PP_B(n)$ or
$\PP_D(n)$.  This will allow us to make concise statements about
$\Pe(n)$ for $\epsilon \in \{0,1\}$.

If $\lambda = [\lambda_1 \ge \cdots \ge \lambda_k]$, we write
$\sigma_j(\lambda)$ for the $i$-th \emph{partial sum} $\sum_{i=1}^j
\lambda_i$.  Recall the standard partial order on partitions: for
$\lambda, \lambda' \in \PP(n)$, we say that $\lambda \le \lambda'$ if
we have $\sigma_j(\lambda) \le \sigma_j(\lambda')$ for all $j$.  In
this case, we say that $\lambda'$ \emph{dominates} $\lambda$.  Recall
also that the closure order on nilpotent orbits coincides with this
order on partitions in the classical groups.  For a partition
$\lambda$, let $\lambda^*$ denote its transpose partition, and let
$\lambda_B$, $\lambda_C$, $\lambda_D$ denote its $B$-, $C$-, and
$D$-collapses respectively, whenever those are defined.  (The
$X$-collapse of $\lambda$ is the unique largest partition $\lambda'$
such that $\lambda' \le \lambda$ and $\lambda' \in \PP_X(n)$; see
\cite{col-mcg:nilp}.)  Suppose $\lambda = [ \lambda_1 \ge \cdots \ge
\lambda_k ] \in \PP(n)$, and assume that $\lambda_k \ne 0$.  We define
the following four operations:
\begin{align*}
\lambda^+ &= [\lambda_1 + 1 \ge \lambda_2 \ge \cdots \ge \lambda_k] & 
\lambda^- &= [\lambda_1 \ge \cdots \ge \lambda_{k-1} \ge \lambda_k - 1] \\
\lambda_+ &= [\lambda_1 \ge \cdots \ge \lambda_k \ge 1] & 
\lambda_- &= \lambda^{*-*}
\end{align*}
(Note that $\lambda_+ = \lambda^{*+*}$ as well.)

Given two partitions $\lambda$ and $\mu$, we can form their
\emph{union} $\lambda \cup \mu$, a partition of $|\lambda| + |\mu|$,
by putting $r_{\lambda \cup \mu}(a) = r_\lambda(a) + r_\mu(a)$ for all
$a$.  We can also take their \emph{join}, defined by
\[
\lambda \vee \mu = (\lambda^* \cup \mu^*)^*.
\]
If one thinks of partitions in terms of Young diagrams, the union
corresponds to combining the rows of the two diagrams, while the join
corresponds to combining their columns.  Finally, if $\lambda =
[\lambda_1 \ge \cdots \ge \lambda_k]$, we define
\[
\chi^+_j(\lambda) = [\lambda_1 \ge \cdots \ge \lambda_j]
\qquad\text{and}\qquad
\chi^-_j(\lambda) = [\lambda_{j+1} \ge \cdots \ge \lambda_k].
\]
Note that $\lambda = \chi^+_j(\lambda) \cup \chi^-_j(\lambda)$ for any
$j$.

Sometimes we will want to restrict the kinds of partitions that we
take unions and joins of, in order to have control over what the union
or join looks like.  Given two partitions $\lambda$ and $\mu$, let $a$
be the smallest part of $\lambda$, and let $b$ be the largest part of
$\mu$.  We say that $\lambda$ is \emph{superior} to $\mu$ if $a \ge
b$.  We say that $\lambda$ is \emph{evenly} (resp.~\emph{oddly})
\emph{superior} to $\mu$ if there is an even (resp.~odd) number $m$
such that $a \ge m \ge b$.

\subsection{Computing with collapses}

The following observations about collapses will be relied upon heavily
when we set about the work of proving the main theorems in
Section~\ref{sect:proof-class}.  If $\lambda$ has $k$ parts, then any
collapse $\lambda_X$ of it must have either $k$ or $k+1$ parts.
Moreover, $B$-partitions necessarily have an odd number of parts, and
$D$-partitions necessarily have an even number, so we can determine
exactly how many parts $\lambda_B$ or $\lambda_D$ must have (of
course, only one of those collapses is defined for any particular
$\lambda$).  Finally, $\lambda_C$ (when it is defined) must have the
same number of parts as $\lambda$, because if it had one more, we
would have introduced a new part equal to $1$, but we cannot create
new odd parts when taking a $C$-collapse.

We will often encounter situations in which we have a partition
written as the union or join of two others, and in which we will want
to express a certain collapse of $\lambda$ in terms of collapses of
the smaller partitions.  The following proposition collects formulas
for twelve kinds of joins, and twelve kinds of unions.  This table of
formulas is certainly sufficient for the calculations in this paper.
The author has not bothered to determine whether any of the
twenty-four could have been omitted.

\begin{lem}\label{lem:cupvee}
Suppose $\lambda = \lambda' \vee \lambda''$.  Let $k$ be the
largest part of $\lambda'$, and let $p = |\lambda'|$.  Suppose in
addition that $\mu = \mu' \cup \mu''$, that $\mu'$ has $k$ parts, and
that $|\mu'| = p$.  Assume that $\lambda'{}^*$ is superior to
$\lambda''{}^*$, and that $\mu'$ is superior to $\mu''$.  The following
table expresses various collapses of $\lambda$ and $\mu$ in terms of
collapses of the smaller partitions.  For any formula containing
$\lambda''{}_B$, we must make the additional assumption that
$\lambda'{}^*$ is oddly superior to $\lambda''{}^*$; for any containing
$\lambda''{}_D$, we assume that $\lambda'{}^*$ is evenly superior to
$\lambda''{}^*$.  Similarly, for any formula containing $\mu'{}^-$ and
$\mu''{}^+$, we  must assume that $\mu'$ is superior to $\mu''{}^+$.
\begin{center}
\begin{tabular}[b]{lcccc}
& \multicolumn{2}{c}{$k$ even}&\multicolumn{2}{c}{$k$ odd}\\
& $p$ even     & $p$ odd     & $p$ even     & $p$ odd     \\
$\lambda_B$: &
  $\lambda'{}^+\!_{B-} \vee \lambda''{}_B$ &
  $\lambda'{}_B        \vee \lambda''{}_D$ &
  $\lambda'{}^+\!_B    \vee \lambda''{}^-\!_C$ &
  $\lambda'{}_B        \vee \lambda''{}_C$  \\
$\lambda_C$: &
  $\lambda'{}_C        \vee \lambda''{}_C$ &
  $\lambda'{}^+\!_C    \vee \lambda''{}^-\!_C$ &
  $\lambda'{}_C        \vee \lambda''{}_D$ &
  $\lambda'{}^+\!_{C-} \vee \lambda''{}_B$  \\
$\lambda_D$: &
  $\lambda'{}_D        \vee \lambda''{}_D$ &
  $\lambda'{}^+\!_{D-} \vee \lambda''{}_B$ &
  $\lambda'{}_D        \vee \lambda''{}_C$ &
  $\lambda'{}^+\!_D    \vee \lambda''{}^-\!_C$ \\
{\tiny\ }\\
$\mu_B$: &
  $\mu'{}_D        \cup \mu''{}_B$ &
  $\mu'{}^-\!_D    \cup \mu''{}^+\!_B$ &
  $\mu'{}^-\!_B    \cup \mu''{}^+\!_D$ &
  $\mu'{}_B        \cup \mu''{}_D$  \\
$\mu_C$: &
  $\mu'{}_C        \cup \mu''{}_C$ &
  $\mu'{}^-\!_C    \cup \mu''{}^+\!_C$ &
  $\mu'{}_C        \cup \mu''{}_C$ &
  $\mu'{}^-\!_C    \cup \mu''{}^+\!_C$  \\
$\mu_D$: &
  $\mu'{}_D        \cup \mu''{}_D$ &
  $\mu'{}^-\!_D    \cup \mu''{}^+\!_D$ &
  $\mu'{}^-\!_B    \cup \mu''{}^+\!_B$ &
  $\mu'{}_B        \cup \mu''{}_B$
\end{tabular}
\end{center}
\end{lem}
\begin{proof}
Once one becomes accustomed to the pattern of producing these
formulas, it is fairly easy to compute all of them.  We will work
through just one: that for $\lambda_B$ when $k$ is odd and $p$ is
even.  For $\lambda_B$ to be defined, $|\lambda|$ must be odd; and
since $p$ is even, $|\lambda''|$ must be odd.  Since $k$ is odd, the
parities of parts of $\lambda''$ are opposite to those of the
corresponding parts of $\lambda$, so taking a $B$-collapse of
$\lambda$ should manifest itself as something like a $C$-collapse of
$\lambda''$.  Since $|\lambda''|$ is odd, if we attempt to take a
$C$-collapse of it, we will be partway through a collapsing operation
when we get to the end of the partition: there will be a leftover
``$1$'' to be added to some odd part, but no remaining odd parts to
receive it.  This ``$1$'' will ``leak'' onto $\lambda'$.  We can
preemptively take care of this leaking $1$ by looking at $\lambda'{}^+$
and $\lambda''{}^-$ instead.  Now, we comfortably take the $C$-collapse
of $\lambda''{}^-$, and the $B$-collapse of $\lambda'{}^+$.  (If
$\lambda''$ has $m$ parts, it may seem that we should have added the
leaking $1$ to the $(m+1)$-th part of $\lambda'$, not its first part,
as is done by writing $\lambda'{}^+$.  But in $\lambda'{}^+$, the first
part is now even, and the remaining parts up to the $m$-th one are all
odd, so in taking a $B$-collapse, that ``$1$'' gets shoved down to at
least the $(m+1)$-th row anyway.)  We thus obtain that $\lambda_B =
\lambda'{}^+\!_B \vee \lambda''{}^-\!_C$.

The only comment we make on other cases is regarding the auxiliary
superiority requirements.  Terms of the form $\lambda''{}_B$ or
$\lambda''{}_D$ may have a different number of parts from $\lambda''$,
so we have to be a lot more careful in considering the interaction
between $\lambda'$ and $\lambda''$.  The easiest thing to do is impose
a condition that the largest part of $\lambda'$ have high enough
multiplicity that we need not worry: that is exactly what the
superiority condition does for us.  Similar considerations result in
the corresponding requirements when we deal with $\mu'{}^-$ and
$\mu''{}^+$. 
\end{proof}

\subsection{Marked partitions}

If $X$ is one of $B$, $C$, or $D$, we define $\tP_X(n)$ to be the set
of pairs of partitions $(\nu,\eta)$, such that:
\begin{enumgen}{1}
\item $\nu \cup \eta \in \PP_X(n)$.
\item Every part of $\nu$ is odd (resp.~even) if $X = B$ or $D$
(resp.~$C$) and has multiplicity $1$.
\item If $X = B$ or $D$, $\nu$ has an even number of parts.
\end{enumgen}
This notation is taken from \cite{sommers:duality}, but we will
typically find another notation far more convenient for our purposes.
We will write elements $(\nu,\eta) \in \tP_X(n)$ as
$\cmp{\lambda}{\nu}$, where $\lambda = \nu \cup \eta$.  In this
notation, we think of elements of $\tP_X(n)$ just as partitions from
$\PP_X(n)$, with the additional data that certain parts ({\it viz.} 
those in $\nu$) have been ``marked.''  Indeed, we will refer to
elements of these sets as \emph{marked partitions}, and we call $\nu$
the \emph{marking partition} and $\lambda$ the \emph{underlying
partition}.  Marked partitions of the form $\cmp{\lambda}{\vn}$ are
called \emph{trivially marked partitions}.  As before, we sometimes
write $\tP_0(n)$ and $\tP_1(n)$ for these sets.

We can attempt to define the union and join operations for marked
partitions, but the constructions we give now may not always yield a
valid marked partition.  This situation will be rectified in the
following subsection, when we introduce ``reduced marked partitions.''
For now, we define the union simply by
\begin{equation}\label{eqn:union-def}
\cmp{\lambda_1}{\nu_1} \cup \cmp{\lambda_2}{\nu_2}
= \cmp{(\lambda_1 \cup \lambda_2)}{(\nu_1 \cup \nu_2)}.
\end{equation}
Next, write $\lambda_1 \vee \lambda_2 = [a_1 \ge \cdots \ge a_k]$.
Suppose $\nu_1 = [n_1 \ge \cdots \ge n_p]$, and $\nu_2 = [m_1 \ge
\cdots \ge m_r]$.  We define
\begin{equation}\label{eqn:join-def}
\cmp{\lambda_1}{\nu_1} \vee \cmp{\lambda_2}{\nu_2}
= \cmp{(\lambda_1 \vee \lambda_2)}{\omega},
\qquad\text{where}\qquad
\omega = [ a_{\height_{\lambda_1}(n_i)} \mid i = 1,\ldots,p ] \cup
[ a_{\height_{\lambda_2}(m_i)} \mid i = 1,\ldots,r ].
\end{equation}
The idea of this definition is that we should preserve the heights of
the marked parts when we take the join.  Quite often, we will
encounter joins of marked partitions in which the largest part of
$\lambda_2$ has very high multiplicity, more than the total number of
parts of $\lambda_1$.  In this special circumstance, understanding the
join of marked partitions is much easier: if $b$ is that largest part
of $\lambda_2$, we obtain
\[
\omega = [ b+n_1 \ge \cdots \ge b+n_p \ge m_1 \ge \cdots \ge m_r ].
\]

For $\lambda \in \Pe(n)$ and $\delta \in \{0,1\}$, let
\[
S_\delta(\lambda) = \{ a \mid \text{$a \not\equiv \epsilon \pmod{2}$
and $r(a) \equiv \delta \pmod{2}$} \}.
\]
We will just write $S_\delta$ when no confusion will result.  For
$\cmp{\lambda}{\nu} \in \tPe(n)$, write
\[
T_\delta(\pmp{\lambda}{\nu}) = T_\delta = \nu \cap S_\delta(\lambda).
\]

\subsection{Parametrizing $\boldsymbol{\N}\!_{\text{\bf o,c}}$ and 
$\boldsymbol{\N}\!_{\text{\bf o,\=c}}$}

A detailed account of the following description of a parametrization
of $\Noc$ and $\Nocc$ can be found in \cite{sommers:b-c}.
Now, $\tP_X(n)$ is close to indexing the set $\Noc$ in type $X$.
Actually, there is a surjective map
\[
\tP_X(n) \to \Noc
\]
which is a bijection in type $B$, but is $2$-to-$1$ over any orbit in
types $C$ and $D$ whose partition has $S_1 \ne \vn$.  There
is, of course, a further projection
\begin{equation}\label{eqn:proj-nocc}
\tP_X(n) \to \Nocc.
\end{equation}
We now describe this projection in some detail.  Given $\lambda$, list
the elements of $S_1$ as $j_l > \cdots > j_1$.  Assume that $l$ is even
in type $C$ by taking $j_1 = 0$ if necessary ($l$ is automatically odd
in type $B$ and even in type $D$).  Now, given $\cmp{\lambda}{\nu}$, let
$T_0^{(m)} = \{ a \in T_0(\pmp{\lambda}{\nu}) \mid j_m < a < j_{m+1}\}$,
and let $T_1^{(m)} = T_1 \cap \{j_m\}$.
Next, we define an equivalence relation $\sim$ on $\tPe(n)$ as
follows: $\cmp{\lambda}{\nu} \sim \cmp{\lambda}{\nu'}$ if
\begin{enumgen}{1}
\item $T_0^{(m)}(\pmp{\lambda}{\nu}) = T_0^{(m)}(\pmp{\lambda}{\nu'})$
whenever $m$ is even. \label{enum:sim1}
\item $\big|T_1^{(m+1)}(\pmp{\lambda}{\nu}) \cup
  T_0^{(m)}(\pmp{\lambda}{\nu}) \cup 
  T_1^{(m)}(\pmp{\lambda}{\nu})\big|   \equiv 
  \big|T_1^{(m+1)}(\pmp{\lambda}{\nu'}) \cup
  T_0^{(m)}(\pmp{\lambda}{\nu'}) \cup
  T_1^{(m)}(\pmp{\lambda}{\nu'})\big| \pmod{2}$
whenever $m$ is odd. \label{enum:sim2}
\end{enumgen}
(In the second of these conditions, we interpret $T_1^{(l+1)}$ as
$\vn$ in type $B$.)  Then, the projection in (\ref{eqn:proj-nocc}) is
precisely the quotient by $\sim$.  We can formulate one particular
equivalence under $\sim$ quite easily, as follows.  If we are working
in type $B$, let $\tilde S_1 = S_1 \setminus \{ j_l \}$, and note that
this set has an even number of elements.

\begin{lem}\label{lem:sim-equiv}
Given a marked partition $\cmp{\lambda}{\nu} \in \tP_X(n)$, define
\[
\nu' = T_0(\pmp{\lambda}{\nu}) \cup \begin{cases}
S_1(\lambda) \setminus T_1(\pmp{\lambda}{\nu}) &
  \text{in types $C$ and $D$,} \\
(\tilde S_1(\lambda) \setminus T_1(\pmp{\lambda}{\nu})) \cup
(T_1(\pmp{\lambda}{\nu}) \cap \{ j_l \}) &
  \text{in type $B$.}
\end{cases}
\]
Then $\cmp{\lambda}{\nu'} \sim \cmp{\lambda}{\nu}$.
\end{lem}
\begin{proof}
It is easy to see that condition~(\ref{enum:sim2}) above is satisfied
when we replace $T_1$ by its complement in $S_1$ in types $C$ and
$D$.  In type $B$, we need to be careful when $m = l$, because there
is no $j_{l+1}$, but the same idea goes through if we take only take
the complement of that portion of $T_1$ which meets $\tilde S_1$, as
in the above formula.
\end{proof}

Consider the set
\[
\tPo_X(n) = \big\{ \cmp{\lambda}{\nu} \in \tP_X(n) \mathbin{\big|}
\text{$T_1^{(m+1)} = T_0^{(m)} = \vn$ whenever $m$ is odd} \big\},
\]
which we call the set of \emph{reduced marked partitions}.  It is easy
to see that the restricted map $\tPo_X(n) \to \Nocc$ is a bijection.
An alternate description of these sets is as follows.  If $\lambda$ is
of type $B$ (resp.~$C$, $D$), let us call a part of $\lambda$
\emph{markable} if it is odd (resp.~even, odd) and has odd
(resp.~even, even) height.  Then we have
\[
\tPo_X(n) = \{ \cmp{\lambda}{\nu} \in \tP_X(n) \mid
\text{$\nu$ consists only of markable parts of $\lambda$} \}.
\]

We will speak of elements of $\tP_X(n)$ as \emph{labels} for elements
of $\Noc$ and $\Nocc$, and of elements of $\tPo_X(n)$ as the
\emph{reduced labels} for elements of $\Nocc$.  Every element of
$\tP_X(n)$ is $\sim$-equivalent to exactly one element of $\tPo_X(n)$.
The process of passing to the reduced label can be described as
follows.  Given $\cmp{\lambda}{\nu}$, we define a new marked partition
$\cmp{\lambda}{\nu'}$, which is characterized as follows: we have
\begin{align*}
T_0^{(m)}(\pmp{\lambda}{\nu'}) &= \begin{cases}
T_0^{(m)}(\pmp{\lambda}{\nu})  & \text{if $m$ is even} \\
\vn            & \text{if $m$ is odd}
\end{cases}\\
\intertext{and}
T_1^{(m)}(\pmp{\lambda}{\nu'}) &= \begin{cases}
\{j_m\}     & \text{if $m$ is odd and
              $\big|T_1^{(m+1)}(\pmp{\lambda}{\nu}) \cup 
              T_0^{(m)}(\pmp{\lambda}{\nu}) \cup
              T_1^{(m)}(\pmp{\lambda}{\nu})\big|$ is odd} \\ 
\vn & \text{otherwise.}
\end{cases}
\end{align*}

There is often a sensible way, given two arbitrary partitions
$\lambda$ and $\nu$, to carry out a ``reduction'' procedure that
generalizes the above one, even when $\cmp{\lambda}{\nu}$ is not a
valid marked partition.  (This goes hand-in-hand with the idea of
generalized height.)  We must first fix one of the types $B$, $C$, or
$D$ as the context in which we are working, but we do not require that
$\lambda$ be a partition of that type.  The only condition we impose
is that when the context type is $B$ or $D$, $\nu$ must have an even
number of parts.  Note that the above definition of ``markable part''
makes sense without any restriction on $\lambda$, and let $m_k >
\cdots > m_1$ be the set of markable parts of $\lambda$.  We define
$\nu'$ by putting
\[
r_{\nu'}(m_i) = \begin{cases}
1 & \text{if $\height_\nu(m_i) - \height_\nu(m_{i+1})$ is odd,} \\
0 & \text{if $\height_\nu(m_i) - \height_\nu(m_{i+1})$ is even,}
\end{cases}
\qquad
\text{$r_{\nu'}(a) = 0$ if $a$ is not a markable part of $\lambda$,}
\]
where, when $i = k$, we interpret $\height_\nu(m_{k+1})$ as $0$.  It
is easy to verify that when $\cmp{\lambda}{\nu}$ is a marked
partition, this generalized procedure coincides with the above one for
passing to a reduced marked partition.

Let us return to the problem that unions and joins of marked
partitions may not yield valid marked partitions.  Typically, we
employ the above procedure to pass from whatever partitions the
formulas~(\ref{eqn:union-def}) and~(\ref{eqn:join-def}) yield to a
reduced marked partition.  Indeed, henceforth, unless explicitly
stated otherwise, all marked partitions are assumed to be reduced, and
if any possibly nonreduced marked partition appears in a formula, we
silently assume that it is to be replaced by an equivalent reduced
one.

\subsection{Duality and special orbits}

We now recall the formulas for $\dLS$, $\dBV$, and $\dS$ in the
classical groups.  
\begin{equation}\label{eqn:dS-def}
\begin{aligned}
\text{Type $B$:}&\qquad& \dLS(\lambda) &= \lambda^*\!_B \qquad&
  \dBV(\lambda) &= \lambda^-\!_C\!^* \qquad&
  \dS(\nu,\eta) &= (\nu \cup \eta^-\!_C)^*\!_C \\
\text{Type $C$:}&\qquad& \dLS(\lambda) &= \lambda^*\!_C \qquad&
  \dBV(\lambda) &= \lambda^+\!_B\!^* \qquad&
  \dS(\nu,\eta) &= (\nu \cup \eta^+\!_B)^*\!_B &\\
\text{Type $D$:}&\qquad& \dLS(\lambda) &= \lambda^*\!_D \qquad&
  \dBV(\lambda) &= \lambda^*\!_D \qquad&
  \dS(\nu,\eta) &= (\nu \cup \eta^*\!_D\!^*)^*\!_D
\end{aligned}
\end{equation}
The formulas for $\dBV$ are obtained by combining the formulas for
$\dLS$ with the following formulas for the order-preserving bijection
between $\Nosp$ and $\Langdual\Nosp$ in types $B$ and $C$:
\begin{alignat*}{2}
\Nosp(B_n) \to \Nosp(C_n) &: \qquad \lambda \mapsto \lambda^-\!_C \\
\Nosp(C_n) \to \Nosp(B_n) &: \qquad \lambda \mapsto \lambda^+\!_B 
\end{alignat*}
In fact, these same formulas can be evaluated on nonspecial partitions
in $\No(B_n)$ and $\No(C_n)$: they then compute the following
composition of maps:
\[
\xymatrix@1{*+{\No} \ar[r]^-{\dLS^2} & *+{\Nosp} \ar[r]^-{\simeq}
&*+{\Langdual\Nosp}}. 
\]
The formulas for $\dS$ are given in \cite{sommers:duality}.  We are
now in a position to revisit the proof of
Proposition~\ref{prop:triv-min}.

\begin{proof}[Proof of Proposition~\ref{prop:triv-min} in the
classical types]
We need to show that $\dS(\orb,1) \ge \dS(\orb,C)$; this should follow
from a quick computation using the above formulas.  We carry it out
now in type $D$.  Starting with $\lambda = \nu \cup
\eta$, it is easy to see that
\begin{subequations}
\begin{align}
                   \eta &\le \eta^*\!_D\!^*  \\
\lambda = \nu \cup \eta &\le \nu \cup \eta^*\!_D\!^* \label{eqn:dS-ordb} \\
              \lambda^* &\ge (\nu \cup \eta^*\!_D\!^*)^* \label{eqn:dS-ordc} \\
          \lambda^*\!_D &\ge (\nu \cup \eta^*\!_D\!^*)^*\!_D. 
                                                  \label{eqn:dS-ordd} 
\end{align}
\end{subequations}
Essentially the same reasoning works in types $B$ and $C$ as well,
although we need to replace~(\ref{eqn:dS-ordb}) above with the
following slightly less trivial inequalities:
\begin{align*}
\lambda^-\!_C &\le \nu \cup \eta^-\!_C, &
\lambda^+\!_B &\le \nu \cup \eta^+\!_B.
\end{align*}
Moreover, in types $B$ and $C$, we need to use the observations that
$\lambda^-\!_C\!^* = \lambda^-\!_C\!^*\!_C$ and $\lambda^+\!_B\!^* =
\lambda^+\!_B\!^*\!_B$, respectively, to pass from~(\ref{eqn:dS-ordc})
to~(\ref{eqn:dS-ordd}).
\end{proof}

We also recall the recipe for computing Sommers' canonical inverse.
if $\lambda$ is of type $B$ (resp.~$C$ or $D$), we let $\pi$ be the
set of even (resp.~odd) parts of $\lambda^*$ with odd multiplicity.
Then the canonical inverse is given by $\cmp{\dBV(\lambda)}{\pi}$,
where we pass to the reduced marked partition if necessary.  (In
\cite{sommers:duality}, Sommers regards the canonical inverse as a map
$\Langdual\No \to \Noc$, so he made no comment about passing a reduced
marked partition, but in the present context, we regard it as a
map $\Langdual\No \to \Nocc$.)

The images of $\dLS$ and $\dBV$ consist precisely of the set of
special orbits, which are labelled by special partitions.  A
characterization of special partitions may be found in
\cite{col-mcg:nilp}.  If $\lambda$ is a $B$- (resp. $C$-, $D$-)
partition, it is special if all its even (resp.~odd, even) parts have
odd (resp.~even, even) height.  Moreover, if $\lambda$ is a special
$B$- (resp.~$C$-) partition, then $\lambda^*$ is also a special $B$-
(resp~$C$-) partition.  If $\lambda$ is a special $D$-partition, then
$\lambda^*$ is a (not necessarily special) $C$-partition.  We conclude
with a lemma about formulas for special partitions.

\begin{lem}\label{lem:special-alt}
The following identities hold:
$\lambda^-\!_C\!^* = \lambda^{*-}\!_C$ 
 for $\lambda \in \PP_B(n)$,
$\lambda^+\!_B\!^* = \lambda^{*+}\!_B$ 
 for $\lambda \in \PP_C(n)$, and
$\lambda^*\!_D\!^* = \lambda^{+-}\!_C$ 
 if either $\lambda \in \PP_D(n)$ or $\lambda^* \in \PP_C(n)$.
\end{lem}
\begin{proof}
The proof establishes all three formulas simultaneously by induction
on the sum of the partition.  One verifies it by direct calculation
for the smallest partitions: $[3]$ and $[1^3]$ in type $B$, $[2]$ and
$[1^2]$ in type $C$, and $[1^2]$ in type $D$.  We work out the
inductive step when $\lambda$ is of type $B$; the others are handled
similarly.  Let $m = \#\lambda$, and let $b$ be the smallest part of
$\lambda$.  Note that $m$ is necessarily odd.  We can write $\lambda =
[b^m] \vee \lambda'$, where $\lambda'$ is a $B$-partition if $b$ is
even, and a $C$-partition if $b$ is odd.  Suppose first that $b$ is
odd.  We have $\lambda^- = [b^m]^- \vee \lambda'$, so $\lambda^-\!_C =
[b^m]^-\!_C \vee \lambda'_D$ by Lemma~\ref{lem:cupvee}.
Now, $[b^m]^-\!_C = [b^{m-1},b-1]$, so we get $\lambda^-\!_C\!^* =
[m^{b-1},m-1] \cup \lambda'_D\!^*$.  Using the inductive hypothesis,
we rewrite this as $[m^b]^- \cup \lambda'{}^{*+-}\!_C$.  Another appeal
to Lemma~\ref{lem:cupvee} lets us conclude that this last expression
is equal to $\lambda^{*-}\!_C$.  The case of $b$ even is handled
similarly, as are types $C$ and $D$. 
\end{proof}

\section{Construction in the classical groups}
\label{sect:constr-class}

We are now ready to define the map $\dA: \Nocc \to \Langdual\Nocc$.
Since we know that we want $\pr1 \circ \dA$ to agree with $\dS$, where
$\pr1 : \Nocc \to \No$ is projection to the first member, the main
difficulty is defining the marking partition on the range.  Given a
reduced marked partition $\cmp{\lambda}{\nu}$, write $\nu = [n_l >
\cdots > n_1]$, and assume that $l$ is even by taking $n_1 = 0$ if
necessary in type $C$.  Define
\begin{equation}\label{eqn:nuhat-def}
\hat\nu = [\height_\lambda(n_1)-1 > \cdots > \height_\lambda(n_l)-1].
\end{equation}
(If we are in type $C$ and $n_1 = 0$, we need to say what
$\height_\lambda(0)$ means.  We want this quantity to be even, since
markable parts are supposed to have even height in type $C$.  We take
it to be the smallest even number larger than $\#\lambda$.)  Next, if
$\cmp{\lambda}{\nu} = (\nu,\eta)$ is a marked partition of type $B$
(resp.~$C$, $D$), we define
\begin{equation}\label{eqn:pi-def}
\pi = \{ \text{even (resp.~odd, odd) parts of $\eta^*$ with odd
multiplicity} \}.
\end{equation}
We regard this set as a partition, each of whose parts has
multiplicity $1$.  We then put
\begin{equation}\label{eqn:dA-def}
\begin{aligned}
\text{\it Type $B$:}&\quad& 
  \cmp{\tau}{\rho} &= \cmp{\nu^*}{\vn} \vee \cmp{\eta^-\!_C\!^*}{\pi}
  &\quad 
  \dA(\pmp{\lambda}{\nu}) &= \cmp{\tau_C}{\hat\nu \cup \rho} \\
\text{\it Type $C$:}&\quad& 
  \cmp{\tau}{\rho} &= \cmp{\nu^*}{\vn} \vee \cmp{\eta^+\!_B\!^*}{\pi}
  &\quad 
  \dA(\pmp{\lambda}{\nu}) &= \cmp{\tau_B}{\hat\nu \cup \rho} \\
\text{\it Type $D$:}&\quad& 
  \cmp{\tau}{\rho} &= \cmp{\nu^*}{\vn} \vee \cmp{\eta^*\!_D}{\pi}
  &\quad 
  \dA(\pmp{\lambda}{\nu}) &= \cmp{\tau_D}{\hat\nu \cup \rho}.
\end{aligned}
\end{equation}
Computing with these formulas by hand can be quite cumbersome,
especially since one must pass to a reduced marked partition several
times.  To help elucidate the matter, we now discuss in detail the
steps involved in computing $\dA$ in type $B$.  First, we find the
partition $\tilde\eta = \eta^-\!_C\!^*$.  Next, $\pi$ may not be a
valid marking partition for $\tilde\eta$, or it may simply not be
reduced: in any case, we reduce $\cmp{\tilde\eta}{\pi}$ to obtain a
reduced marked partition $\cmp{\tilde\eta}{\pi'}$.  Now, $\tau$ is
simply the join $\nu^* \vee \tilde\eta$.  To compute $\rho$, we must
refer to the description of joins of marked partitions in
Section~\ref{sect:partitions}: we first take $\rho'$ to be the set of
parts of $\tau$ whose heights are the same as the heights of the parts
of $\pi$ in $\tilde\eta$.  Then, $\cmp{\tau}{\rho}$ is the reduced
marked partition obtained by reducing $\cmp{\tau}{\rho'}$.  Finally,
let $\tilde\tau$ be the collapse $\tau_C$.  Again, $\hat\nu \cup \rho$
may not be a valid reduced marking partition for $\tilde\tau$, but
reducing $\cmp{\tilde\tau}{\hat\nu \cup \rho}$ will give us some
$\cmp{\tilde\tau}{\sigma}$.  This is the final answer:
$\dA(\pmp{\lambda}{\nu}) = \cmp{\tilde\tau}{\sigma}$.  The following
example illustrates this procedure.

\begin{exam}
Consider the orbit $\orb$ labelled by $[7,5,4^2,3,2^2,1^2]$ in type
$B_{14}$, or $\son(29)$.  This partition has three markable parts:
$7$, $3$, and $1$.  Therefore, $\Ab(\orb) \simeq (\Z/2\Z)^2$; the four
possible marking partitions are $\vn$, $[3,1]$, $[7,3]$, and $[7,1]$.
Let us consider the conjugacy class corresponding to $[3,1]$.  Writing
$\cmp{[7,5,4^2,3,2^2,1^2]}{[3,1]}$ as a pair, we have $(\nu,\eta) =
([3,1], [7,5,4^2,2^2,1])$.  We compute $\eta^-\!_C = [6^2,4^2,2^2]$,
which is self-dual: $\eta^-\!_C\!^* = [6^2,4^2,2^2]$.  We have $\nu^*
= [2,1^2]$, so $\tau = [8,7,5,4,2^2]$.  Finally, $\tau_C =
[8,6^2,4,2^2]$.

To compute the marking partition, we have $\eta^* = [7,6,4^2,2,1^2]$,
so $\pi = [6,2]$.  Both parts of $\pi$ are markable in
$\eta^-\!_C\!^*$, so $\cmp{\eta^-\!_C\!^*}{\pi}$ is already reduced.
Taking the join with $\nu^*$ yields $\cmp{[8,7,5,4,2^2]}{[7,2]}$,
which becomes $\cmp{[8,7,5,4,2^2]}{[4,2]}$ when we reduce it.
Finally, $\hat\nu = [8,4]$, so for the final answer, we take the
reduced marked partition corresponding to
$\cmp{[8,6^2,4,2^2]}{[8,4^2,2]}$, arriving at
\[
\dA(\pmp{[7,5,4^2,3,2^2,1^2]}{[3,1]}) = \cmp{[8,6^2,4,2^2]}{[4,2]}.
\]
\end{exam}

In this and the following section, we will be ensconced in many
laborious computations with the above formulas.  Most of the results
must actually be proved thrice, once in each of types $B$, $C$, and
$D$; but we will usually only write out the full details in type $B$,
and just make cursory remarks about the nature of the calculations in
the other types.

We can establish the following two properties of $\dA$ immediately
from the definition.  Once again, $\pr1: \Nocc \to \No$ is the
obvious projection map.  Let us also recall Sommers' canonical
inverse, which was mentioned near the end of
Section~\ref{sect:formal}.  This is a certain right inverse to $\dS$;
see \cite{sommers:duality} for its construction.

\begin{prop}\label{prop:dS-resp}
We have that $\pr1 \circ \dA$ agrees with $\dS$.
\end{prop}
\begin{proof}
In type $B$, the underlying partition of $\dA(\pmp{\lambda}{\nu})$ is
$\tau_C = (\nu^* \vee \eta^-\!_C\!^*)_C = (\nu \cup
\eta^-\!_C)^*\!_C$, which is precisely the formula for $\dS(\nu,\eta)$
in type $B$.  Types $C$ and $D$ are equally easy to handle.
\end{proof}

\begin{prop}\label{prop:dS-inv-class}
Given an orbit $\orb$ labelled by a partition $\lambda$, the conjugacy
class labelled by $\dA(\pmp{\lambda}{\vn})$ coincides with
Sommers' canonical inverse for $\orb$.
\end{prop}
\begin{proof}
When $\nu = \vn$ and $\eta = \lambda$, the formula for
$\cmp{\tau}{\rho}$ in~(\ref{eqn:dA-def}) agrees with Sommers' recipe
for the canonical inverse.  At this stage, $\tau = \dBV(\lambda)$ is
already a $C$-, $B$-, or $D$-partition (in types $B$, $C$, and $D$,
respectively), so the additional collapse of $\tau$ in the formula for
$\dA$ does nothing.  We also have $\hat\nu = \vn$, so
$\dA(\pmp{\lambda}{\nu}) = \cmp{\tau}{\pi}$.
\end{proof}

Before we can set about proving that the above map is, in fact, an
extended duality map as defined in Section~\ref{sect:intro}, we need
to develop some techniques for manipulating marked partitions.  The
formulas we have so far are too opaque to be tackled in their raw form
when we want to prove things about them.  We spend the rest of the
section showing how to break down a marked partition into ``blocks,''
and how to compute $\dA$ piecemeal on the individual blocks.

\begin{lem}\label{lem:square-blocks}
Suppose that $\cmp{\lambda}{\nu} = \cmp{[a^l]}{\vn} \vee
\cmp{\lambda'}{\nu'}$, where $\lambda$ has $l$ parts.
\begin{enumgen}{1}
\item\label{enum:square-blocks-BD}
If $\cmp{\lambda}{\nu}$ is of type $B$ or $D$ and $a$ is even, then
$\dA(\pmp{\lambda}{\nu}) = \cmp{[l^a]}{\vn} \cup
\dA(\pmp{\lambda'}{\nu'})$.  In this case, $\cmp{\lambda'}{\nu'}$ is
of the same type as $\cmp{\lambda}{\nu}$.
\item\label{enum:square-blocks-C}
If $\cmp{\lambda}{\nu}$ is of type $C$, $a$ is odd, and $l$ is even,
then $\dA(\pmp{\lambda}{\nu}) = \cmp{[l+1,l^{a-1}]}{\vn} \cup
\dA(\pmp{\lambda'}{\nu'})$.   Here, $\cmp{\lambda'}{\nu'}$ is of type
$D$.
\end{enumgen}
In both cases, we also have $\dA^2(\pmp{\lambda}{\nu}) =
\cmp{[a^l]}{\vn} \vee \dA^2(\pmp{\lambda'}{\nu'})$.
\end{lem}
\begin{proof}
Let us write $\cmp{\lambda}{\nu} = (\nu,\eta)$ and
$\cmp{\lambda'}{\nu'} = (\nu',\eta')$, and let us refer back to the
formulas for $\dS$.  If $\eta$ has $n$ parts, then
\[
\nu = [a^{l-n}] \vee \nu'
\qquad\text{and}\qquad
\eta = [a^n] \vee \eta'.
\]
We will prove part~(\ref{enum:square-blocks-BD}) when
$\cmp{\lambda}{\nu}$ is of type $B$.  The type-$D$ case of
part~(\ref{enum:square-blocks-BD}), as well as
part~(\ref{enum:square-blocks-C}) and the statement for $\dA^2$, are
handled similarly.  Now, $\eta^-$ may be given by either $[a^n] \vee
\eta'{}^-$ or $[a^{n-1},a-1] \vee \eta'$, depending on whether $\eta'$
has $n$ parts or fewer than $n$ parts.  We compute $\eta^-\!_C$ with
the appropriate formula from Lemma~\ref{lem:cupvee}, and see
that in either case, we get $[a^n] \vee \eta'{}^-\!_C$ (possibly using
the fact that $[a^{n-1},a-1]^+\!_C = [a^n]$).  Therefore,
\[
\nu \cup \eta^-\!_C = ([a^{l-n}] \vee \nu') \cup ([a^n] \vee \eta'{}^-\!_C)
                    = [a^l] \vee (\nu' \cup \eta'{}^-\!_C).
\]
Taking the transpose of both sides, we get
\[
(\nu \cup \eta^-\!_C)^* = [l^a] \cup (\nu' \cup \eta'{}^-\!_C)^*.
\]
Now, we use another formula from Lemma~\ref{lem:cupvee} to
compute the $C$-collapse of this expression.  We obtain
\begin{equation}\label{eqn:square-blocks-part}
\dS(\pmp{\lambda}{\nu}) = (\nu \cup \eta^-\!_C)_C = [l^a] \cup (\nu' \cup
\eta'{}^-\!_C)_C = [l^a] \cup \dS(\pmp{\lambda'}{\nu'}).
\end{equation}
We now need to compute the marking partition.  If $\pi$ is defined
from $\eta$ according to~(\ref{eqn:pi-def}), and $\pi'$ is defined
analagously from $\eta'$, we evidently have $\pi = \pi'$.  It is then
easy to work through the formulas of~(\ref{eqn:dA-def}) and see that
$\rho = \rho'$ as well, and finally that $\dA(\pmp{\lambda}{\nu})$ and
$\dA(\pmp{\lambda'}{\nu'})$ have the same marking partition.
\end{proof}

\begin{lem}\label{lem:B-blocks}
Given $\cmp{\lambda}{\nu} \in \tP_B(n)$, suppose that
$\cmp{\lambda}{\nu} = \cmp{\lambda_1}{\nu_1} \cup
\cmp{\lambda_2}{\nu_2}$, with $\cmp{\lambda_1}{\nu_1} \in \tP_B(m)$
and $\cmp{\lambda_2}{\nu_2} \in \tP_D(n-m)$.  Suppose furthermore that
$\lambda_1$ is evenly superior to $\lambda_2$.  Then 
$
\dA(\pmp{\lambda}{\nu}) = \dA(\pmp{\lambda_1}{\nu_1}) \vee
\dA(\pmp{\lambda_2}{\nu_2})
$.
(Note that the term $\dA(\pmp{\lambda_2}{\nu_2})$ is to be computed in
type $D$.)
\end{lem}
\begin{proof}
Write $\cmp{\lambda}{\nu} = (\nu,\eta)$, $\cmp{\lambda_1}{\nu_1} =
(\nu_1,\eta_1)$, and $\cmp{\lambda_2}{\nu_2} = (\nu_2,\eta_2)$ in the
notation of pairs, and let $\pi$, $\pi_1$, and $\pi_2$ be the
corresponding partitions as defined in~(\ref{eqn:pi-def}).  Note that
$\eta_1$ has an odd number of parts, and $|\eta_1|$ is odd, both
because $\eta_1$ is a $B$-partition.  We consult
Lemma~\ref{lem:cupvee}, starting with $\eta^- = \eta_1 \cup \eta_2^-$,
and find that $\eta^-\!_C = \eta_1^-{}_C \cup \eta_2^{-+}\!_C$; then,
Lemma~\ref{lem:special-alt} tells us that we actually have $\eta^-\!_C
= \eta_1^-{}_C \cup \eta_2^*{}_D\!^*$.  Because $\eta_1$ has an odd
number of parts, even parts of $\eta^*$ correspond to odd parts of
$\eta_2^*$, so $\pi_1$ and $\pi_2$ are related to $\pi$ as in the
following equation:
\[
\cmp{\eta^-\!_C\!^*}{\pi} = \cmp{\eta_1^-{}_C\!^*}{\pi_1} \vee
\cmp{\eta_2^*{}_D}{\pi_2}.
\]
Since $\nu = \nu_1 \cup \nu_2$, we have $\nu^* = \nu_1^* \vee \nu_2^*$;
it follows directly that
\begin{equation}\label{eqn:tau-blocks}
\cmp{\nu^*}{\vn} \vee \cmp{\eta^-\!_C\!^*}{\pi} =
  (\cmp{\nu_1^*}{\vn} \vee \cmp{\eta_1^-{}_C\!^*}{\pi_1}) \vee
  (\cmp{\nu_2^*}{\vn} \vee \cmp{\eta_2^*{}_D}{\pi_2}).
\end{equation}
Write this equation, following~(\ref{eqn:dA-def}), as
$\cmp{\tau}{\rho} = \cmp{\tau_1}{\rho_1} \vee \cmp{\tau_2}{\rho_2}$.
Now, $\tau_1$ has an odd number of columns, because $\#\lambda_1$ is
odd.  (One might worry that it could have fewer columns due to the
``$^-$'' operation, if the smallest part of $\eta_1$ were $1$, but
that is not possible since $\lambda_1$ is superior to $\lambda_2$.)
Now we make use of the hypothesis of even superiority: another appeal
to Lemma~\ref{lem:cupvee} tells us exactly that $\tau_C =
\tau_1{}_C \vee \tau_2{}_D$.  At this point, we have established that
\begin{equation}\label{eqn:dS-blocks}
\dS(\pmp{\lambda}{\nu}) = \dS(\pmp{\lambda_1}{\nu_1}) \vee
\dS(\pmp{\lambda_2}{\nu_2}).
\end{equation}

Now consider the marking partition.  If $n$ is a part of $\nu_2$, we
have $\height_\lambda(n) = \height_{\lambda_2}(n) + \#\lambda_1$.  We
therefore have $\hat\nu = \hat\nu_1 \cup
([(\#\lambda_1)^{\#\hat\nu_2}] \vee \hat\nu_2)$.  Combining this
description with~(\ref{eqn:dS-blocks}) and~(\ref{eqn:tau-blocks}), we
find
\[
\cmp{\tau_C}{\hat\nu \cup \rho} = \cmp{\tau_1{}_C}{\hat\nu_1 \cup \rho_1}
\vee \cmp{\tau_2{}_D}{\hat\nu_2 \cup \rho_2},
\]
as desired.
\end{proof}

\noindent
Entirely analagous arguments establish the three cases of the
following lemma.

\begin{lem}\label{lem:CD-blocks}
Let $\cmp{\lambda}{\nu}$ be a marked partition, and suppose that
$\cmp{\lambda}{\nu} = \cmp{\lambda_1}{\nu_1} \cup \cmp{\lambda_2}{\nu_2}$.
\begin{enumgen}{1}
\item \label{enum:CD-blocks-C}
If $\cmp{\lambda}{\nu} \in \tP_C(n)$, let us also suppose that
$\cmp{\lambda_1}{\nu_1} \in \tP_C(m)$, that $\cmp{\lambda_2}{\nu_2} \in
\tP_C(n-m)$, that $\lambda_1$ and $\nu_1$ have an even number of parts,
and that $\lambda_1$ is oddly superior to $\lambda_2$.  Then
$
\dA(\pmp{\lambda}{\nu}) = \dA(\pmp{\lambda_1}{\nu_1})_- \vee
\dA(\pmp{\lambda_2}{\nu_2})
$.
Here, $\dA(\pmp{\lambda_1}{\nu_1})_-$ is to be understood as applying
the $_-$ operation to the underlying partition.  The largest part of
the underlying partition is odd, so the marking partition is
unaffected.
\item \label{enum:CD-blocks-D}
If $\cmp{\lambda}{\nu} \in \tP_D(n)$, let us also suppose that
$\cmp{\lambda_1}{\nu_1} \in \tP_D(m)$, that $\cmp{\lambda_2}{\nu_2} \in
\tP_D(n-m)$, and that $\lambda_1$ is evenly superior to
$\lambda_2$. Then
$
\dA(\pmp{\lambda}{\nu}) = \dA(\pmp{\lambda_1}{\nu_1}) \vee
\dA(\pmp{\lambda_2}{\nu_2})
$.
\item \label{enum:CD-blocks-C-aux}
In the context of case~(\ref{enum:CD-blocks-C}), let us further
suppose that $\lambda$, $\nu$, $\lambda_2$, and $\nu_2$ have even
numbers of parts.  Then
$
\dA(\pmp{\lambda}{\nu})_- = \dA(\pmp{\lambda_1}{\nu_1})_- \vee
\dA(\pmp{\lambda_2}{\nu_2})_-
$. \qed
\end{enumgen}
\end{lem}

\noindent
Case~(\ref{enum:CD-blocks-C-aux}) of this lemma may seem bizarre, but
we will arrive at a use for it shortly.

Returning to the context of Lemma~\ref{lem:B-blocks}, let $m$ be the
even number arising in the definition of ``evenly superior'' for
$\lambda_1$ and $\lambda_2$, and let $l = \#\lambda_1$.  (Note that
$l$ is odd.)  There is marked $B$-partition $\cmp{\lambda'_1}{\nu'_1}$
such that $\cmp{\lambda_1}{\nu_1} = \cmp{[m^l]}{\vn} \vee
\cmp{\lambda'_1}{\nu'_1}$.  Using Lemma~\ref{lem:square-blocks}, we can
write
\begin{align}
\dA(\pmp{\lambda}{\nu}) &= (\cmp{[l^m]}{\vn} \cup
\dA(\pmp{\lambda'_1}{\nu'_1})) \vee \dA(\pmp{\lambda_2}{\nu_2}) \notag
\\
&= (\cmp{[l^m]}{\vn} \vee \dA(\pmp{\lambda_2}{\nu_2}))
\cup \dA(\pmp{\lambda'_1}{\nu'_1}) , \label{eqn:dA2-initial}
\end{align}
where we have made use of the fact that the largest part of
$\lambda_2$ is at most $m$, so the underlying partition of
$\dA(\pmp{\lambda_2}{\nu_2})$ has at most $m$ parts.  We can apply
part~(\ref{enum:square-blocks-C}) of Lemma~\ref{lem:square-blocks} to
the first term and write
\[
\dA(\cmp{[l^m]}{\vn} \vee \dA(\pmp{\lambda_2}{\nu_2})) = 
\cmp{[m+1,m^{l-1}]}{\vn} \cup \dA^2(\pmp{\lambda_2}{\nu_2}).
\]
Now, the union in~(\ref{eqn:dA2-initial}) is exactly of the form
demanded by part~(\ref{enum:CD-blocks-C}) of
Lemma~\ref{lem:CD-blocks}, so we can apply that statement here.
\begin{align*}
\dA^2(\pmp{\lambda}{\nu}) &= \dA(\cmp{[l^m]}{\vn} \vee
\dA(\pmp{\lambda_2}{\nu_2}))_- \vee \dA^2(\pmp{\lambda'_1}{\nu'_1}) \\
&= (\cmp{[m^l]}{\vn} \cup \dA^2(\pmp{\lambda_2}{\nu_2})) \vee
\dA^2(\pmp{\lambda'_1}{\nu'_1}) \\
&= (\cmp{[m^l]}{\vn} \vee \dA^2(\pmp{\lambda'_1}{\nu'_1})) \cup
\dA^2(\pmp{\lambda_2}{\nu_2}).
\end{align*}
By one final application of Lemma~\ref{lem:square-blocks}, we obtain
the following result for type $B$.  Similar calculations establish it
in types $C$ and $D$.

\begin{lem}\label{lem:d2-blocks}
Let $\cmp{\lambda}{\nu} = \cmp{\lambda_1}{\nu_1} \cup
\cmp{\lambda_2}{\nu_2}$ be a decomposition as in
Lemma~\ref{lem:B-blocks} or~\ref{lem:CD-blocks}.  Then
$\dA^2(\pmp{\lambda}{\nu}) = \dA^2(\pmp{\lambda_1}{\nu_1}) \cup
\dA^2(\pmp{\lambda_2}{\nu_2})$. \qed
\end{lem}

Now, we can use Lemma~\ref{lem:CD-blocks}(\ref{enum:CD-blocks-D})
iteratively to split up a marked $D$-partition into smaller and
smaller pieces.  We can also do the same in type $B$, if we use
Lemma~\ref{lem:B-blocks} and
Lemma~\ref{lem:CD-blocks}(\ref{enum:CD-blocks-D}) in combination.
Similarly, parts~(\ref{enum:CD-blocks-C})
and~(\ref{enum:CD-blocks-C-aux}) of Lemma~\ref{lem:CD-blocks} taken
together let us split up $C$-partitions into smaller and smaller
pieces.  The following definition captures the precise nature of the
permitted decompositions.

\begin{defn}\label{defn:blocks}
Suppose we have $\cmp{\lambda}{\nu} = \cmp{\lambda_1}{\nu_1} \cup
\cdots \cup \cmp{\lambda_k}{\nu_k}$.  Such a decomposition is called a
\emph{division into blocks} of $\cmp{\lambda}{\nu}$, and each
$\cmp{\lambda_i}{\nu^i}$ is called a \emph{block}, under the
circumstances described below.

If $\cmp{\lambda}{\nu} \in \tP_B(n)$, we require that
$\cmp{\lambda_1}{\nu_1} \in \tP_B(k_1)$ and that
$\cmp{\lambda_i}{\nu_i} \in \tP_D(k_i)$ for $i > 1$.  Furthermore,
$\lambda_i$ must be evenly superior to $\lambda_{i+1}$ for $i = 1,
\ldots, k-1$.

If $\cmp{\lambda}{\nu} \in \tP_C(n)$, we require that
$\cmp{\lambda_i}{\nu_i} \in \tP_C(k_i)$ for all $i$, and that
$\lambda_i$ and $\nu_i$ have an even number of parts for $i = 1,
\ldots, k-1$.  Furthermore, $\lambda_i$ must be oddly superior to
$\lambda_{i+1}$ for $i = 1, \ldots, k-1$.

If $\cmp{\lambda}{\nu} \in \tP_D(n)$, we require that
$\cmp{\lambda_i}{\nu_i} \in \tP_D(m)$ for all $i$.  Furthermore,
$\lambda_i$ must be evenly superior to $\lambda_{i+1}$ for $i = 1,
\ldots, k-1$.
\end{defn}

\noindent
We now combine Lemmas~\ref{lem:B-blocks}, \ref{lem:CD-blocks},
and~\ref{lem:d2-blocks} to obtain the following concise statement.

\begin{prop}\label{prop:blocks}
Let $\cmp{\lambda}{\nu} = \cmp{\lambda_1}{\nu_1} \cup \cdots \cup
\cmp{\lambda_k}{\nu_k}$ be a division into blocks. Then, $\dA$ can be
computed as follows: 
\begin{align*}
\text{\it Type $B$:}& &
\dA(\pmp{\lambda}{\nu}) &= \dA(\pmp{\lambda_1}{\nu_1}) \vee \cdots \vee
\dA(\pmp{\lambda_k}{\nu_k}) \\
\text{\it Type $C$:}& &
\dA(\pmp{\lambda}{\nu}) &= \dA(\pmp{\lambda_1}{\nu_1})_- \vee \cdots \vee
\dA(\pmp{\lambda_{k-1}}{\nu_{k-1}})_- \vee \dA(\pmp{\lambda_k}{\nu_k}) \\
\text{\it Type $D$:}& &
\dA(\pmp{\lambda}{\nu}) &= \dA(\pmp{\lambda_1}{\nu_1}) \vee \cdots \vee
\dA(\pmp{\lambda_k}{\nu_k})\qed
\end{align*}
Moreover, in all types, $\dA^2(\pmp{\lambda}{\nu}) =
\dA^2(\pmp{\lambda_1}{\nu_1}) \cup \cdots \cup
\dA^2(\pmp{\lambda_k}{\nu_k})$. \qed 
\end{prop}

The motivation for developing the idea of divisions into blocks is our
hope that we can cut up arbitrary marked partitions into blocks that
are very simple in some sense, and that such blocks will be easy to
work with when we set about the task of proving the main theorems.
We now state precisely the sort of blocks we hope to obtain.

\begin{defn}\label{defn:basic}
A \emph{basic block} of type $B$ (resp.~$C$, $D$) is a marked
partition $\cmp{\lambda}{\nu}$ such that $\nu$ has one or two parts,
say $[n_2]$ or $[n_2 > n_1]$, such that $n_1$ (if it exists) is the
smallest part of $\lambda$, and such that $n_2$ is the largest part of
odd (resp.~even, even) height in $\lambda$.  The circumstance of $\nu$
having only one part can occur only in type $C$; in this case, we
often regard $\nu$ as having two parts by putting $n_1 = 0$.  A basic
block is called \emph{ultrabasic} if it meets the additional condition
that $n_1 \le 1$.
\end{defn}

\begin{prop}
Any marked partition has a division into blocks such that each block
is either a trivially marked partition or a basic block.
\end{prop}
\begin{proof}
This is easily seen by induction on the number of parts of the
underlying partition.  Given a marked partition $\cmp{\lambda}{\nu}$
of type $B$, let $a$ be the first part of odd height.  If $a$ is even
(and therefore unmarkable) or odd and unmarked, we put $\lambda_1 =
\chi^+_{\height(a)}(\lambda)$ and $\lambda_2 =
\chi^-_{\height(a)}(\lambda)$.  Then, $\cmp{\lambda}{\nu} =
\cmp{\lambda_1}{\vn} \cup \cmp{\lambda_2}{\nu}$ is a division
into blocks in which the first term is trivially marked, and in the
second, $\lambda_2$ has fewer parts than $\lambda$.

If $a$ is a marked part, let $b$ be the second marked part.  This time
we take $\lambda_1 = \chi^+_{\height(b)}(\lambda)$ and $\lambda_2 =
\chi^-_{\height(b)}(\lambda)$.  This time, $\cmp{\lambda}{\nu} =
\cmp{\lambda_1}{[a>b]} \cup \cmp{\lambda_2}{\nu \setminus \{a,b\}}$ is
a division into blocks whose first term is a basic block.  Similar
arguments work in types $C$ and $D$.
\end{proof}

Henceforth, all our arguments regarding properties of $\dA$ will
address only basic and trivially marked blocks.

\begin{prop}\label{prop:basic}
If $\cmp{\lambda}{\nu}$ is a basic block, then
$\dS(\pmp{\lambda}{\nu})$ can be computed by the following simplified
formulas: 
\begin{alignat*}{2}
\text{Type $B$:}&\quad& \dS(\pmp{\lambda}{\nu}) &= \lambda^{-*}\!_C \\
\text{Type $C$:}&\quad& \dS(\pmp{\lambda}{\nu}) &= \lambda^{+*}\!_B \\
\text{Type $D$:}&\quad& \dS(\pmp{\lambda}{\nu}) &= \lambda^{+-*}\!_D
\end{alignat*}
\end{prop}
\begin{proof}
If $\cmp{\lambda}{\nu}$ is a basic block, write $\nu = [n_2 > n_1]$.
Let us assume for the time being that $\cmp{\lambda}{\nu}$ is an
ultrabasic block.  This will make our calculations less cumbersome.
We will obtain a formula; then, at the end of the proof, we use
Lemma~\ref{lem:square-blocks} to see that the same formula holds for
general basic blocks.

Suppose we are working in type $B$, so $n_1 = 1$.  Let $h_i =
\height_\lambda(n_i)$ for $i = 1,2$.  Thus $h_1$ is the total number
of parts of $\lambda$.  Let $\mu_1 = \chi^+_{h_2-1}(\lambda)$ and
$\mu_2 = \chi^-_{h_2-1}(\lambda)$.  Note that $\mu_1$ has only parts
of even height, and that $r_{\mu_2}(n_2) = 1$.  Let $\mu_2'$ be the
partition gotten from $\mu_2$ by decreasing the multiplicities of
$n_2$ and $n_1$ each by $1$.  Since $n_1 = 1$, we have $\mu_2 = [n_2]
\cup \mu_2'{}_+$.  Writing $\cmp{\lambda}{\nu}$ as a pair
$(\nu,\eta)$, we have $\eta = \mu_1 \cup \mu_2'$, and $\eta^- = \mu_1
\cup \mu_2'{}^-$.  Using Lemma~\ref{lem:cupvee}, we get $\eta^-\!_C =
\mu_1{}_C \cup \mu_2'{}^-\!_C$, but since all parts of $\mu_1$ have even
height, they all have even multiplicity, so $\mu_1{}_C = \mu_1$:
\begin{equation}\label{eqn:basic-part1}
\eta^-\!_C = \mu_1 \cup \mu_2'{}^-\!_C.
\end{equation}
Next, again using that $n_1 = 1$, we have $\nu \cup \eta^-\!_C =
\mu_1 \cup [n_2] \cup \mu_2'{}^-\!_C{}_+$, or
\begin{equation}\label{eqn:basic-part2}
(\nu \cup \eta^-\!_C)^* = \mu_1^* \vee [n_2]^* \vee
\mu_2'{}^-\!_C\!^{*+}.
\end{equation}
We use Lemma~\ref{lem:cupvee} to get $(\nu \cup
\eta^-\!_C)^*\!_C = \mu_1^*{}_C \vee ([n_2]^* \vee
\mu_2'{}^-\!_C\!^{*+})_C$.  Since $\mu_1$ only has parts of even height,
$\mu_1^*$ only has even parts, so the $C$-collapse does nothing.
Using Lemma~\ref{lem:cupvee} yet again, we find that the second
term is equal to $[n_2]^* \vee \mu_2'{}^-\!_C\!^{*+}\!_B$.  Now $\mu_2'$
is a $B$-partition, so Spaltenstein's formulas give us that
$\mu_2'{}^-\!_C\!^{*+}\!_B = \mu_2'{}^*\!_B$.  Finally, using the fact
that $[n_2]^* \vee \mu_2'{}^* = \mu_2^{-*}$, we obtain $[n_2]^* \vee
\mu_2'{}^*\!_B = \mu_2^{-*}\!_C$:
\begin{equation}\label{eqn:basic}
(\nu \cup \eta^-\!_C)^*\!_C = \mu_1^* \vee \mu_2^{-*}\!_C.
\end{equation}
Now, we know $\lambda^{-*} = \mu_1^* \vee \mu_2^{-*}$, and
Lemma~\ref{lem:cupvee} would tell us that $\lambda^{-*}\!_C =
\mu_1^*{}_C \vee \mu_2^{-*}\!_C$.  But since $\mu_1$ only has parts of
even height, $\mu_1^*$ only has even parts, and the $C$-collapse does
nothing do it.  Thus~(\ref{eqn:basic}) is given by $\lambda^{-*}\!_C$,
as desired.

We do not give the details in types $C$ and $D$, but as an aid to
those who wish to work them out, we list the analogues
of (\ref{eqn:basic-part1}), (\ref{eqn:basic-part2}),
and~(\ref{eqn:basic}) here.
\begin{align*}
\begin{gathered}
\text{\it Type $C$}\\
[n_2] \cup \eta^+\!_B = \mu_1 \cup \mu_2^+{}_{B-} \\
(\nu \cup \eta^+\!_B)^* = \mu_1^* \vee \mu_2^+{}_B\!^{*-}  \\
(\nu \cup \eta^+\!_B)^*\!_B = \mu_1^{+*} \vee \mu_2^*{}_C
\end{gathered}
&&
\begin{gathered}
\text{\it Type $D$} \\
\eta^*\!_D\!^* = \mu_1^+ \cup \mu_2'{}^{*-}\!_C\!^* \\
(\nu \cup \eta^*\!_D\!^*)^* = \mu_1^{+*} \vee [n_2]^* \vee
\mu_2'{}^{*-}\!_C\!^+ \\
(\nu \cup \eta^*\!_D\!^*)^*\!_D = \mu_1^{+*} \vee \mu_2^{-*}\!_C
\end{gathered}
\end{align*}
In type $C$, it turns out to be more convenient not to work with
$\mu_2'$.  In type $D$, we need to make use of the identity
$\eta^*\!_D\!^* = \eta^{+-}\!_C$.  With these points in mind, the
proofs are straightforward.
\end{proof}

\section{Proofs of the main theorems in the classical groups}
\label{sect:proof-class}

In this section, we establish the main theorems of the paper for the
classical groups.  Theorem~\ref{thm:po} is relatively easy: we prove
it first, and we make use of it from time to time as we go about
proving Theorem~\ref{thm:d}.  The proof of the latter is broken up
into a number of steps and occupies most of the section.  The steps
may look familiar: we end up proving that $\dA$ has many of the
properties established in Section~\ref{sect:formal} before we show
that it is actually an extended duality map.

\subsection{The partial order in the classical groups}

The strategy for the proof of the theorem below is quite simple: we
just attempt the raw computation of the two values of $\dS$, using the
techniques from the previous section.  Those techniques make it
straightforward to find a difference in the answers, starting with a
difference in the original marking partitions.

\begin{thm}\label{thm:po-class}
Let $C, C' \subset A(\orb)$ be two conjugacy classes associated to the
same orbit.  Then, in the classical groups, $\dS(\orb,C) =
\dS(\orb,C')$ if and only if $C$ and $C'$ have the same image in
$\Ab(\orb)$.  As a consequence, the partial order (\ref{eqn:po}) is
well-defined.
\end{thm}
\begin{proof}
We need to prove that if $C$ and $C'$ are two different conjugacy
classes in $\Ab(\orb)$, then $\dS(\orb,C) \ne \dS(\orb,C')$.  Suppose
that these conjugacy classes are labelled by $\cmp{\lambda}{\nu}$ and
$\cmp{\lambda}{\nu'}$, respectively.  Let $a$ be the largest part of
$\lambda$ that appears in only one of $\nu$ and $\nu'$.  Therefore,
$a$ has (generalized) heights of opposite parity in $\nu$ and $\nu'$;
assume it has even height in $\nu$.  That means that we can break
$\cmp{\lambda}{\nu}$ up into blocks $\cmp{\lambda_1}{\nu_1} \cup
\cmp{\lambda_2}{\nu_2}$, where the smallest part of $\lambda_1$ is $a$.
(Note that because $a$ is markable, this is a legitimate division into
blocks in whatever type we are working in.)  But in $\cmp{\lambda}{\nu'}$,
there is some basic block $\cmp{\zeta}{\omega}$, $\omega = [w_2 > w_1]$,
such that $w_2 \ge a > w_1$.  We build a division into blocks around
this basic block, writing $\cmp{\lambda}{\nu'} =
\cmp{\lambda'_1}{\nu'_1} \cup \cmp{\zeta}{\omega} \cup
\cmp{\lambda'_2}{\nu'_2}$.  Finally, let $h = \height_\lambda(a)$.

Let $\mu = \dS(\pmp{\lambda}{\nu})$, and $\mu' =
\dS(\pmp{\lambda}{\nu'})$.  Using Proposition~\ref{prop:blocks} just
to compute $\dS$, we have
\[
\mu^* = \dS(\pmp{\lambda_1}{\nu_1})^* \cup \dS(\pmp{\lambda_2}{\nu_2})^*.
\]
We see that $\sigma_h(\mu^*) = |\dS(\pmp{\lambda_1}{\nu_1})^*| = |\lambda_1|
= \sigma_h(\lambda)$ in types $C$ or $D$, and $\sigma_h(\mu^*) =
|\lambda_1| - 1 = \sigma_h(\lambda) - 1$ in type $B$.  (This comes
from just counting the ``$^+$'' and ``$^-$'' operations that are done
in computing $\dS$ in each type.)

We now analyze $\mu'$.  Write $\zeta = \zeta' \cup \zeta''$, where
$\zeta' = \chi^+_{\height(a)}(\zeta)$ and $\zeta'' =
\chi^-_{\height(a)}(\zeta)$.  Suppose we are working in type $B$; and
suppose further that $\cmp{\lambda'_1}{\nu'_1}$ is nontrivial, so that
$\cmp{\zeta}{\omega}$ is of type $D$.  (Definition~\ref{defn:blocks}
says basic $B$-blocks can only occur at the beginning of a marked
$B$-partition.)  Then Proposition~\ref{prop:basic} says
$\dS(\pmp{\zeta}{\omega}) = \zeta^{+-*}\!_D$.  We have $\zeta^{+-*} =
\zeta'{}^{+*} \vee \zeta''{}^{-*}$, so by Lemma~\ref{lem:cupvee},
$\zeta^{+-*}\!_D = \zeta'{}^{+*+}\!_{D-} \vee \zeta''{}^{-*}\!_B$, so
\[
\mu'{}^* = \dS(\pmp{\lambda'_1}{\nu'_1}) \cup \zeta'{}^{+*+}\!_{D-} \cup
\zeta''{}^{-*}\!_B \cup \dA(\pmp{\lambda'_2}{\nu'_2}).
\]
We see that $\sigma_h(\mu'{}^*) = |\dS(\pmp{\lambda'_1}{\nu'_1})| +
|\zeta'{}^{+*+}\!_{D-}| = (|\lambda'_1| - 1) + (|\zeta'| + 1) =
|\lambda'_1| + |\zeta'| = \sigma_h(\lambda)$.  A nearly identical
argument establishes that $\sigma_h(\mu'{}^*) = \sigma_h(\lambda)$ when
$\cmp{\lambda'_1}{\nu'_1}$ is trivial and $\cmp{\zeta}{\omega}$ is of
type $B$.

Similar computations show that in types $C$ and $D$, we get
$\sigma_h(\mu'{}^*) = \sigma_h(\lambda) + 1$.  Thus, in every case, we
get $\sigma_h(\mu'{}^*) = \sigma_h(\mu^*) + 1$, so $\mu \ne
\mu'$, as desired.
\end{proof}

\subsection{Special marked partitions}

The remainder of the section is devoted to establishing
Theorem~\ref{thm:d} for the classical groups.  We begin our attack on
it by attempting to characterize the marked partitions that occur in
the image of $\dA$.  In this subsection, we define the set
$\tPsp_X(n)$ of special marked partitions, and show that the image of
$\dA$ is contained with this set.  Of course, we still have to prove
various properties of $\dA$ before we can know that this terminology
coincides with the idea of ``special'' that we introduced in
Section~\ref{sect:formal}.

\begin{defn}\label{defn:special}
Let $\cmp{\lambda}{\nu}$ be a reduced marked partition, with $\nu =
[n_l > \cdots > n_1]$.  Assume that $l$ is even, if necessary by
taking $n_1 = 0$ in type $C$.  In type $B$ (resp.~$C$, $D$),
$\cmp{\lambda}{\nu}$ is called \emph{special} if there are no even
(resp.~odd, even) parts of odd (resp.~even, even) height between
$n_{2i}$ and $n_{2i-1}$ for $i = 1, \ldots, l/2$; that is, if there
are no even (resp.~odd, even) parts of odd (resp.~even, even) height
whose (generalized) height in $\nu$ is odd.  The set of special marked
partitions in $\tP_X(n)$ is denoted $\tPsp_X(n)$.
\end{defn}

Note that any trivially marked partition is special by this
definition, as we expect from Proposition~\ref{prop:orbit-special}.
On the other hand, if $\lambda$ is a special partition, a nontrivially
marked partition $\cmp{\lambda}{\nu}$ may be either special or
nonspecial.

\begin{lem}\label{lem:triv-special}
We have that $\dA(\pmp{\lambda}{\vn})$ is a special marked partition
for any $\lambda$.
\end{lem}
\begin{proof}
Let us consider the situation in type $B$.  Recalling
Lemma~\ref{lem:special-alt}, we can write $\dS(\pmp{\lambda}{\vn}) =
\dBV(\lambda) = \lambda^{*-}\!_C$.  Let $\pi$ be the list of even
parts with odd multiplicity in $\lambda^*$.  In $\lambda^{*-}\!_C$,
some parts of $\pi$ have odd multiplicity, and others have even
multiplicity.  According to Lemma~\ref{lem:sim-equiv}, we could
replace $\pi$ by the set $\pi'$ obtained by taking the complementary
set of even parts with odd multiplicity, together with the same set of
even parts with even multiplicity as $\pi$.  That is, $\pi'$ is the
list of even parts that have even multiplicity in $\lambda^{*-}\!_C$
and odd multiplicity in $\lambda^*$, or odd multiplicity in
$\lambda^{*-}\!_C$ and even multiplicity in $\lambda^*$.  We have
$\dS(\pmp{\lambda}{\nu}) = \cmp{\lambda^{*-}\!_C}{\pi} =
\cmp{\lambda^{*-}\!_C}{\pi'}$; we work with $\pi'$ for the rest of
this proof.

What happens as we pass from $\lambda^{*-}$ to $\lambda^{*-}\!_C$?  We
have to make a change in the partition every time we encounter an odd
part with odd multiplicity.  There are an even number of odd parts
with odd multiplicity; we consider them in pairs.  Indeed, suppose
\begin{equation}\label{eqn:pre-collapse}
a_1^{k_1},\,a_2^{k_2},\, \ldots,\, a_l^{k_l}
\end{equation}
is a list of consecutive odd parts of $\lambda^{*-}$, with $a_1$ and
$a_l$ being, say, the largest two odd parts with odd multiplicity.
(We are not requiring that $a_i$ and $a_{i+1}$ be consecutive in
$\lambda^{*-}$, but merely that any parts between them be even.)
We have assumed that the multiplicities $k_1$ and $k_l$ are odd, while
$k_2, \ldots, k_{l-1}$ are even.  Then, the $C$-collapse replaces the
above parts by the following ones:
\begin{equation}\label{eqn:post-collapse}
a_1^{k_1-1},a_1-1,\, a_2+1,a_2^{k_2-2},a_2-1,\, \ldots,\,
a_l+1,a_l^{k_l-1}.
\end{equation}
On each such pair of odd parts with odd multiplicity, the $C$-collapse
follows the pattern of the change from~(\ref{eqn:pre-collapse})
to~(\ref{eqn:post-collapse}); we just investigate what happens on one
instance of the pattern.  Listing the even parts that have changed
multiplicities, we obtain
\begin{equation}\label{eqn:triv-special-mark}
\pi' = [a_1-1,a_2+1,a_2-1, \ldots, a_l+1].
\end{equation}
Note that because $\lambda^{*-}$ has no odd parts between $a_{i-1}$
and $a_i$, there are no odd parts between $a_{i-1}-1$ and $a_i+1$ in
$\lambda^{*-}\!_C$.  

Now, the key observation here is that $\lambda^{*-}\!_C$ is a
special $C$-partition (it equals $\dBV(\lambda)$).  Odd parts in
special $C$-partitions have even height and even multiplicity, so any
part immediately greater than an odd part also has even height.  In
particular, each $a_i+1$ has even height, and is therefore markable.
When we pass to the reduced label, for each $i$, we either retain
both parts $a_{i-1}-1$ and $a_i+1$ (if $a_{i-1}-1$ has even height),
or eliminate both of them (if $a_{i-1}-1$ has odd height).  Since
there are no odd parts between these two parts for any $i$, we have
a special marked partition.
\end{proof}

\begin{prop}\label{prop:im-special}
We have that $\dA(\pmp{\lambda}{\nu})$ is a special marked partition
for any $\cmp{\lambda}{\nu}$.
\end{prop}
\begin{proof}
The previous lemma establishes this fact for trivially marked blocks,
so now we need only consider basic blocks.  This is easy to deduce
from the formulas given in Proposition~\ref{prop:basic}; we work it
out in type $B$ now.  Let $\cmp{\lambda}{\nu}$ be a type-$B$ basic block
with $\nu = [n_2 > n_1]$, and let $h_1 = \height_\lambda(n_1)$: we
have that $h_1$ is odd.  We can write $\lambda = [1^{h_1}] \vee
\lambda'$, where $\lambda'$ is a $C$-partition (since $\lambda$ is a
$B$-partition).  Then, $\lambda^{-*}\!_C = ([h_1-1] \cup
\lambda'{}^*)_C$.  (Note that $[h_1-1]$ is probably not superior to
$\lambda'{}^*$.)  Now, $h_1-1$ is even, so it is unaffected by the
$C$-collapse: $\lambda^{-*}\!_C = [h_1-1] \cup \lambda'{}^*\!_C$.  Since
$\lambda'$ is a $C$-partition, $\lambda'{}^*\!_C$ is a special
$C$-partition, in which all odd parts have even height.  We claim that
the part $[h_1 -1]$ ``pushes them down'' so that they have odd height.
Indeed, the only part of of $\lambda'{}^*\!_C$ larger than $h_1 - 1$ is
$h_1$.  That is, to be sure, an odd part with even height (which is
equal to its multiplicity, $n_2-1$), but all other parts of
$\lambda'{}^*\!_C$ have their heights increased by $1$ when we pass to
$\lambda^{-*}\!_C$.  But since $h_1$ is the largest part of
$\lambda^{-*}\!_C$, it obviously cannot have odd generalized height in
the marking partition.  All other odd parts have odd height, so
$\dA(\pmp{\lambda}{\nu})$ is special.
\end{proof}

\subsection{Involutivity}

Next, we undertake the task of showing that $\dA$ is an involution on
the set of special marked partitions.  We do this in several stages,
beginning just with the trivial conjugacy class on special orbits,
then working up to the trivial conjugacy class for all orbits, and
finally to the full special set.  

\begin{lem}\label{lem:spec-orb-inv}
If $\lambda$ is a special partition, then $\dA(\pmp{\lambda}{\vn}) =
\cmp{\dBV(\lambda)}{\vn}$.
\end{lem}
\begin{proof}
If $\lambda$ is a special $B$- (resp.~$C$-, $D$-) partition, then
$\lambda^*$ is a $B$- (resp.~$C$-, $C$-) partition, so all its even
(resp.~odd, odd) parts have even multiplicity.  Therefore, the
partition $\pi$ defined in~(\ref{eqn:pi-def}) is trivial.  It follows
that $\rho$ in~(\ref{eqn:dA-def}) is trivial as well, as is the
marking partition of $\dA(\pmp{\lambda}{\nu})$.  The underlying
partition is then given by $\dS(\pmp{\lambda}{\nu}) = \dBV(\lambda)$.
\end{proof}

\begin{lem}\label{lem:orb-inv}
For any partition $\lambda$, we have that $\dA^2(\pmp{\lambda}{\vn})
= \cmp{\lambda}{\vn}$.
\end{lem}
\begin{proof}
We know that $\dA(\pmp{\lambda}{\vn})$ is Sommers' canonical
inverse for $\lambda$, so that $\dS(\dA(\pmp{\lambda}{\vn})) =
\lambda$.  We therefore have $\dA^2(\pmp{\lambda}{\vn}) =
\cmp{\lambda}{\nu}$ for some marking partition $\nu$.  We only need to
show that $\nu = \vn$.  To do this, we use
Proposition~\ref{prop:blocks} to decompose $\lambda$ into pieces as
simple as possible.  Call a partition of the form $[a^l]$, with $l$
even, a \emph{rectangle}, and call a partition of the form
\begin{equation}\label{eqn:staircase}
[a_1^{k_1}, a_2^{k_2}, \ldots, a_m^{k_m}],
\end{equation}
with $k_1$ and $k_m$ odd and $k_2, \ldots, k_{m-1}$ even, a
\emph{staircase}.  Any $D$-partition can be written as a union of
rectangles and staircases.  Let us define a \emph{partial staircase}
to be a staircase from which either $a_1^{k_1}$ or $a_m^{k_m}$ is
omitted: a \emph{lower partial staircase} in the former case, and an
\emph{upper partial staircase} in the latter.  In seeking a division
into blocks, we can write any $B$-partition as a lower partial
staircase followed by some number of rectangles and staircases, and
any $C$-partition as a union of rectangles and staircases, possibly
followed by an upper partial staircase.  The proof of this lemma is
accomplished by proving it separately for each of these kinds of
blocks.

In many cases, showing that $\nu$ is trivial is
easy because $\lambda$ just does not have any possible marking
partitions.  Rectangles and staircases have only one part of even
height, while upper partial staircases have none; and lower partial
staircases have only one part of odd height.  The lemma follows in
completely when $\dA(\pmp{\lambda}{\vn})$ is of type $B$ or $D$ (where
marking partitions must have an 
even number of parts), and for upper partial staircases in type $C$.
For rectangles in type $C$, the statement is a consequence of
Lemma~\ref{lem:square-blocks}.  The only remaining case is that of
staircases in type $C$, which we treat now.

Let $\lambda$ be a $C$-partition of the form~(\ref{eqn:staircase}),
and let $\mu = \dS(\pmp{\lambda}{\vn})$.  We have $\mu
=\lambda^+\!_B\!^*$; moreover, we claim that $\dA(\pmp{\lambda}{\nu})
= \cmp{\mu}{\vn}$.  Note that in $\lambda$, $a_1$ and $a_m$ must be
even.  All the parts $a_1, \ldots, a_{m-1}$ have odd height; only
$a_m$ has even height.  Therefore, $\lambda^*$ has only one even part,
its largest one, and the multiplicity of that part is $a_m$, which is
even.  Therefore, $\pi$ as defined in~(\ref{eqn:pi-def}) is trivial,
so $\dA(\pmp{\lambda}{\nu})$ is trivially marked.

We could just trudge ahead and compute $\dA(\pmp{\mu}{\vn})$ directly,
but instead, we use the following trick.  We have already observed
that the proposition holds for $B$-partitions, so
$\dA^2(\pmp{\mu}{\vn}) = \cmp{\mu}{\vn}$.  That is,
$\dA(\pmp{\lambda}{\nu}) = \cmp{\mu}{\vn}$.  But we also have
$\dA(\pmp{\lambda}{\vn}) = \cmp{\mu}{\vn}$, so if we had $\nu \ne
\vn$, that would contradict Theorem~\ref{thm:po-class}.  Therefore,
$\nu = \vn$, and the proposition holds for $C$-staircases.
\end{proof}

\begin{lem}\label{lem:spec-ubasic}
Let $\cmp{\lambda}{\nu}$ be a special ultrabasic block, with $\nu =
[n_2 > n_1]$.  If $\cmp{\lambda}{\nu}$ is of type $B$ (resp.~$C$,
$D$), let $m$ be the largest part of even (resp.~odd, even) height in
$\dS(\pmp{\lambda}{\nu})$.  (When $\cmp{\lambda}{\nu}$ is of type $B$,
we put $m = 0$ if $\dS(\pmp{\lambda}{\nu})$ has no parts of even
height.)  If $m > 1$, then $\dA(\pmp{\lambda}{\nu})$ is again a
special ultrabasic block, given by the following formulas:
\begin{alignat*}{2}
\text{Type $B$:}&\quad& \dA(\pmp{\lambda}{\nu}) &= 
  \cmp{\lambda^{-*}}{[m]}, \\
\text{Type $C$:}&\quad& \dA(\pmp{\lambda}{\nu}) &= 
  \cmp{\lambda^{+*}}{[m,1]}, \\
\text{Type $D$:}&\quad& \dA(\pmp{\lambda}{\nu}) &= 
  \cmp{\lambda^{+-*}}{[m,1]}.
\end{alignat*}
If $m \le 1$, then $\dA(\pmp{\lambda}{\nu})$ is a trivially marked
partition, whose underlying partition is as given above.
\end{lem}
\begin{proof}
To prove this, we must dive into the details of the proof of
Proposition~\ref{prop:im-special}.  Suppose $\cmp{\lambda}{\nu}$ is of
type $B$; recall that we wrote $\lambda = [1^{h_1}] \vee \lambda'$,
where $\lambda'$ is a $C$-partition.  If $\cmp{\lambda}{\nu}$ is
special, then all its even parts must have even heights.  (This is
true for even parts smaller than $n_2$ by the definition of
``special,'' and for even parts larger than $n_2$ by the definition of
``basic block.'')  This means that in $\lambda'$, all odd parts have
even heights; {\it i.e.}, $\lambda'$ is a special $C$-partition.
Therefore, $\lambda'{}^*\!_C = \lambda'{}^*$, and $\lambda^{-*}\!_C =
[h_1-1] \cup \lambda'{}^* = \lambda^{-*}$.  (Here we have used the fact
that $n_1 = 1$, so $h_1$ does not occur as a part of $\lambda'{}^*$.)

Now, $\lambda^* = [h_1] \cup \lambda'{}^*$ does not have any odd parts
with odd multiplicity other than $h_1$, because $\lambda'{}^*$ is a
$C$-partition.  We have $\lambda^* = \eta^{*+} \vee [1^{n_2}]$, so in
$\eta^*$, there are no even parts of height less than or equal to
$n_2$ that have odd multiplicity.  Thus $\pi$ as defined
in~(\ref{eqn:pi-def}) only has parts of height greater than $n_2$.
Now, back in $\lambda$, all parts larger than $n_2$ have even height,
so in $\eta^*$ or $\lambda^*$, all parts whose height is greater than
$n_2$ must be even.  Let $b$ be the largest part of $\lambda^*$ that
has even height greater than $n_2$; in other words, $b$ is the largest
markable part smaller than $h_2$.  (If there are no markable parts
smaller than $h_2$, take $b = 0$.)  There is an odd number of even
parts that are greater than or equal to $b$ and smaller than $h_2$,
but an even number of even parts smaller than $b$ and greater than or
equal to any smaller markable part.  It follows that when we pass to
the reduced marked partition to compute $\cmp{\tau}{\rho}$
in~(\ref{eqn:dA-def}), we get $\rho = [b]$.

Finally, we look at $\dA(\pmp{\lambda}{\nu}) =
\cmp{\lambda^{-*}}{\hat\nu \cup [b]}$.  Since $b$ is the largest
markable part smaller than $h_2$, it is clear that we can replace
$[h_1-1 > h_2 -1 \ge b]$ by $[h_1-1]$ without changing the reduced
marked partition to which it is equivalent.  Now, let $m$ be the
largest part of even height in $\lambda^{-*}$, or, if there are no
parts of even height, take $m = 0$.  If $m \ne 0$, suppose its height
is $k$.  That means that in $\lambda$, $k$ is an even part of height
$m$.  Since $\cmp{\lambda}{\nu}$ is special, $k$ must have even
height, so $m$ is necessarily even.  Therefore, $m$ is markable in
$\lambda^{-*}$.  In the case of either $m \ne 0$ or $m = 0$, then, we
see that $\cmp{\lambda^{-*}}{[m]}$ is the reduced marked partition
equivalent to $\cmp{\lambda^{-*}}{[h_1-1]}$.
\end{proof}

\begin{prop}\label{prop:spec-inv}
If $\cmp{\lambda}{\nu}$ is a special marked partition, then
$\dA^2(\pmp{\lambda}{\nu}) = \cmp{\lambda}{\nu}$.
\end{prop}
\begin{proof}
Lemma~\ref{lem:orb-inv} established this fact for trivially marked
partitions, so now we only need to consider basic blocks.  Indeed, we
actually restrict ourselves to ultrabasic blocks, since we can then
use Lemma~\ref{lem:square-blocks} to pass up to the result for
arbitrary basic blocks.  Let $\cmp{\lambda}{\nu}$ be a special
ultrabasic block, and let $\dA(\pmp{\lambda}{\nu}) = \cmp{\mu}{\xi}$.
According to Lemma~\ref{lem:spec-ubasic}, there are two cases to
consider: either $\cmp{\mu}{\xi}$ is trivially marked, or it is again
a special ultrabasic block.

First, suppose it is trivially marked.  In each type, we can directly
compute $\dS(\pmp{\mu}{\vn}) = \dBV(\mu)$: in type $B$, for
example, we have $\mu = \lambda^{-*}$, and $\dBV(\mu) =
\mu^+\!_B\!^*\!_B = \lambda^*\!_B\!^*\!_B$.  (To get $\mu^+ =
\lambda^*$, we had to use the fact that $n_1 = 1$.)  Moreover,
according to Lemma~\ref{lem:spec-ubasic}, the fact that
$\dA(\pmp{\lambda}{\nu})$ is trivially marked means that $\mu$ has no
parts of even height, which in turn means that $\mu^* = \lambda^-$ has
no even parts.  Again using that $n_1 = 1$, it follows that $\lambda$
has no even parts, and is therefore automatically a special
$B$-partition.  We deduce that $\lambda^*\!_B\!^*\!_B = \lambda$.  Is
it possible that $\dA(\pmp{\mu}{\vn}) = \cmp{\lambda}{\nu'}$ for some
$\nu' \ne \nu$?  Let us again use the trick from the end of the proof
of Lemma~\ref{lem:orb-inv}.  We know from Lemma~\ref{lem:orb-inv} that
$\dA^2(\pmp{\mu}{\vn}) = \cmp{\mu}{\vn}$, but having
$\dA(\pmp{\lambda}{\nu'}) = \cmp{\mu}{\vn} = \dA(\pmp{\lambda}{\nu})$
for $\nu' \ne \nu$ would contradict Theorem~\ref{thm:po-class}.  Thus,
$\dA^2(\pmp{\lambda}{\nu}) = \cmp{\lambda}{\nu}$.

Now, suppose instead that $\xi \ne \vn$.  This time,
$\cmp{\mu}{\xi}$ is itself a special ultrabasic block, so we can use
the formulas of Lemma~\ref{lem:spec-ubasic} twice in a row to
establish the result.  For instance, starting in type $B$, we have
$\mu = \lambda^{-*}$, so $\dA(\pmp{\mu}{\xi}) =
\cmp{\mu^{+*}}{[p,1]}$, where $p$ is the largest part of odd height in
$\mu^{+*}$.  But $\mu^{+*} = \lambda$, as argued in the previous
paragraph, and $n_2$ is the largest part of odd height in $\lambda$.
In this case as well, we find that $\dA^2(\pmp{\lambda}{\nu}) =
\cmp{\lambda}{\nu}$.
\end{proof}

\subsection{Specialization}

The third step is to define a map for passing from a given marked
partition to special one that is larger than it in the partial order.
After we show that that this coincides with $\dA^2$, we will be in a
position prove that $\dA$ is a weak extended duality map, {\it i.e.},
it satisfies axioms~(\ref{ax:ord-rev})--(\ref{ax:dS-resp}).  We begin
with a map which we call the \emph{partial specialization map} $s:
\tPe(n) \to \tPe(n)$, defined as follows.  If $\cmp{\lambda}{\nu}$ is
a nonspecial marked partition of type $B$ (resp.~$C$, $D$), let $a$ be
the smallest even (resp.~odd, even) part of odd (resp.~even, even)
height in $\lambda$ and odd height in $\nu$.  (Of course, no such $a$
exists for a special marked partition).  The part $a$ must have even
multiplicity, say $l$.  Let $\lambda'$ be the partition gotten from
$\lambda$ by deleting all $l$ copies of $a$.  We put
\[
s(\pmp{\lambda}{\nu}) = \begin{cases}
\cmp{\lambda}{\nu}             & 
  \text{if $\cmp{\lambda}{\nu}$ is special,} \\
\cmp{(\lambda' \cup [a+1,a^{l-2},a-1])}{\nu} &
  \text{if $\cmp{\lambda}{\nu}$ is nonspecial.}
\end{cases}
\]
Of course, we may have to pass to the reduced marked partition from
the above formula, if it happens that $a+1$ was a markable part of
$\lambda$ and was, in fact, marked.

It is clear that for a nonspecial marked partition, the map $s$
decreases the total number of even (resp.~odd, even) parts of odd
(resp.~even, even) height in the underlying partition and odd height
in the marking partition.  By induction on that quantity, we obtain
the following result.

\begin{lem}
Given a marked partition $\cmp{\lambda}{\nu}$, there is some
nonnegative integer $N$ such that $s^N(\pmp{\lambda}{\nu})$ is
special. \qed
\end{lem}

We now define the specialization map $e: \tPe(n) \to \tPspe(n)$ as
\[
e(\pmp{\lambda}{\nu}) = s^N(\pmp{\lambda}{\nu}), 
\]
where $N$ is taken large enough that the right-hand side is special.
Note that since the map $s$ fixes special marked partitions, there is
no ambiguity in the above definition arising from the particular
choice of $N$.

\begin{prop}\label{prop:dA-dAs}
We have that $\dA(\pmp{\lambda}{\nu}) = \dA(s(\pmp{\lambda}{\nu}))$
for any marked partition $\cmp{\lambda}{\nu}$.
\end{prop}
\begin{proof}
We begin by proving that $\dS(\pmp{\lambda}{\nu}) =
\dS(s(\pmp{\lambda}{\nu}))$.  Let us assume that $\cmp{\lambda}{\nu}$
is a nonspecial basic block of type $B$, let $a$ and $l$ be as in the
definition of $s$, and let $h = \height_\lambda(a)$.  Let
$\cmp{\zeta}{\omega} = s(\pmp{\lambda}{\nu}) = (\omega,\kappa)$, and
let $b$ be the next smaller part of $\lambda$ after $a$.  Since $a$ is
the smallest even part with odd height, $b$ must either be odd or have
even height.  But if $b$ has even height, it must have odd
multiplicity (since $a$ has odd height), so $b$ is necessarily odd in
all cases.  Suppose $\nu = [n_2 > n_1]$.  If $n_2 > a+1$, then $\omega
= \nu$; otherwise, $\omega = [m > n_1]$, where $m$ is the largest odd
part of $\lambda$ that is smaller than $a$ and has odd height.

Now, $h$ is an odd part with even height $a$ in $\lambda^*$, and $h-l$
(also odd) is the next smaller part after $h$.  Write $\lambda^{-*} =
\tau_1 \cup [h^{a-b},h-l] \cup \tau_2$, where $\tau_1 =
\chi^+_b(\lambda^{-*})$ and $\tau_2 =
\chi^-_{a-1}(\lambda^{-*})$.  Using Lemma~\ref{lem:cupvee}, it
is easy to check that
\[
\dS(\pmp{\lambda}{\nu}) = \lambda^{-*}\!_C = \tau_1{}_C \cup
[h^{a-b-1},h-1,h-l+1] \cup \tau_2{}_C.
\]

We get $\zeta$ from $\lambda$ by replacing $[a^l]$ by $[a+1, a^{l-2},
a-1]$.  Then $\zeta^*$ looks like $\lambda^*$, except that the portion
of the form $[h^{a-b},h-l]$ has been changed to
$[h^{a-b-1},h-1,h-l+1]$.  In the case that $\omega = \nu$, we just
compare with the above computation to see that $\lambda^{-*}\!_C =
\zeta^{-*}\!_C$; {\it i.e.}, $\dS(\pmp{\lambda}{\nu}) =
\dS(s(\pmp{\lambda}{\nu}))$.  But even if $\omega \ne \nu$, we recall
that it is not necessary to pass to reduced marked partition when
computing $\dS$ (which is, after all, defined as a map $\Noc \to
\Langdual\No$), so we can simply replace $\omega$ by $\nu$ and apply
the above argument anyway.

It remains to verify that $\dA(\pmp{\lambda}{\nu})$ and
$\dA(s(\pmp{\lambda}{\nu}))$ produce the same marking partition.  This
is straightforward but extremely tedious.  The proof consists of
writing down the various intermediate marked partitions occurring
in~(\ref{eqn:dA-def}), while scrupulously remembering to pass to a
reduced marked partition whenever possible.  The cases of $\omega =
\nu$ and $\omega \ne \nu$ must be considered separately; the former is
slightly easier.  We omit the details.
\end{proof}

\begin{cor}
We have that $\dA \circ e = \dA$, and that $\dA^2 = e$.
\end{cor}
\begin{proof}
The first statement is an immediate consequence of the preceding
proposition, since, by induction, we have $\dA = \dA \circ s^N$ for
all $N \ge 0$.  Then, on the one hand, we can apply $\dA$ to both
sides again to obtain $\dA^2 \circ e = \dA^2$; but on the other hand,
we know by Proposition~\ref{prop:spec-inv} that $\dA^2$ is the
identity on special marked partitions, and the image of $e$ consists
of special marked partitions, so $\dA^2 \circ e = e$.  Thus $\dA^2 =
e$.
\end{proof}

\begin{prop}\label{prop:ord-pres}
We have that $\dA^2(\pmp{\lambda}{\nu}) \ge \cmp{\lambda}{\nu}$ for
any marked partition $\cmp{\lambda}{\nu}$.
\end{prop}
\begin{proof}
It is easy to see, by construction, that the underlying partition of
$s(\cmp{\lambda}{\nu})$ dominates $\lambda$.  Combining this with
Proposition~\ref{prop:dA-dAs}, we see that $s(\pmp{\lambda}{\nu}) \ge
\cmp{\lambda}{\nu}$.  It follows that $e(\pmp{\lambda}{\nu}) \ge
\cmp{\lambda}{\nu}$; {\it i.e.} that $\dA^2(\pmp{\lambda}{\nu}) \ge
\cmp{\lambda}{\nu}$, as desired.
\end{proof}

\begin{lem}\label{lem:pre-ord-rev}
Suppose that $\cmp{\lambda}{\nu} \le \cmp{\lambda'}{\nu'}$, and that
$\cmp{\lambda'}{\nu'}$ is special.  Then $s(\pmp{\lambda}{\nu}) \le
\cmp{\lambda'}{\nu'}$ as well.
\end{lem}
\begin{proof}
The argument used to prove this statement is similar in flavor to the
argument we gave for Theorem~\ref{thm:po-class}.  Assume that
$\cmp{\lambda}{\nu}$ is nonspecial.  Let $\mu =
\dS(\pmp{\lambda}{\nu})$ and $\mu' = \dS(\pmp{\lambda'}{\nu'})$.  We
have that $\lambda \le \lambda'$ and $\mu \ge \mu'$.  Since
$\dS(s(\pmp{\lambda}{\nu})) = \dS(\pmp{\lambda}{\nu})$, all we have to
prove is that the underlying partition of $s(\pmp{\lambda}{\nu})$ is
smaller than $\lambda'$.  Let $\zeta$ denote the underlying partition
of $s(\pmp{\lambda}{\nu})$.  Let $a$, $l$, and $h$ be as in the proof
of Proposition~\ref{prop:dA-dAs}.  A brief consideration of how
$\zeta$ is formed reveals the following relationship between $\zeta$
and $\lambda$:
\begin{alignat}{2}
\sigma_{h-l+i}(\zeta) &= \sigma_{h-l+i}(\lambda) + 1
&\qquad&\text{for $i = 1, \ldots, l-1$} \notag \\
    \sigma_{k}(\zeta) &= \sigma_{k}(\lambda) 
&\qquad&\text{for $k \ne h-l+1, \ldots, h-1$} \notag \\
\intertext{%
We know $\sigma_k(\lambda) \le \sigma_k(\lambda')$ for all $k$, but to
establish $\zeta \le \lambda'$, we need to prove the following
stronger statements:}
\sigma_{h-l+i}(\lambda) + 1 &\le \sigma_{h-l+i}(\lambda')
&\qquad&\text{for $i = 1, \ldots, l-1$} \label{eqn:sigma-compare}
\end{alignat}
Let us assume that the above fails for some $i$; we shall derive a
contradiction.  Suppose, in particular, that it fails for $i = j$.
This means that $\sigma_{h-l+j}(\lambda) = \sigma_{h-l+j}(\lambda')$.
Let $b_1 \ge \cdots \ge b_l$ be the $(h-l+1)$-th, \dots, $h$-th parts
of $\lambda'$, respectively.  We have
\begin{equation}\label{eqn:sigma-check-j}
\sigma_{h-l+j+1}(\lambda) = \sigma_{h-l+j}(\lambda) + a \le
\sigma_{h-l+j+1}(\lambda') = \sigma_{h-l+j}(\lambda') + b_{j+1},
\end{equation}
so $a \le b_{j+1}$.  We also have
\begin{equation}\label{eqn:sigma-check-0}
\sigma_{h-l+j-1}(\lambda) = \sigma_{h-l+j}(\lambda) - a \le
\sigma_{h-l+j-1}(\lambda') = \sigma_{h-l+j}(\lambda') - b_j,
\end{equation}
which implies $a \ge b_j$.  Since $b_j \ge b_{j+1} \ge a$, we conclude
that $a = b_j = b_{j+1}$.  But then~(\ref{eqn:sigma-check-j}) says
that $\sigma_{h-l+j+1}(\lambda) = \sigma_{h-l+j+1}(\lambda')$,
so~(\ref{eqn:sigma-compare}) fails for $i = j+1$ as well.  If $i_0$ is
the smallest value of $i$ for which~(\ref{eqn:sigma-compare}) fails,
we see by induction that it fails for $i = i_0+1, \ldots, l$ as well.
Furthermore, $b_{i_0} = \cdots = b_l = a$.

We claim, moreover, that $i_0 = 1$; {\it i.e.}
that~(\ref{eqn:sigma-compare}) fails for all $i$.  If not, the
inequality~(\ref{eqn:sigma-check-0}) can be strengthened using the
fact that the $(i_0-1)$-th inequality in~(\ref{eqn:sigma-compare})
holds:
\begin{equation}\label{eqn:sigma-check-1}
\sigma_{h-l+i_0-1}(\lambda) + 1 = \sigma_{h-l+i_0}(\lambda) - a + 1 \le
\sigma_{h-l+i_0-1}(\lambda') = \sigma_{h-l+i_0}(\lambda') - b_{i_0}.
\end{equation}
We deduce that $a - 1 \ge b_{i_0}$.  Since $b_{i_0} \ge b_{i_0+1} \ge
a$, we obtain $a - 1 \ge a$, a contradiction.  We thus have
$\sigma_{h-l+i}(\lambda) = \sigma_{h-l+i}(\lambda')$ and $b_i = a$ for
$i = 1, \ldots, l$.  Additionally, (\ref{eqn:sigma-check-0}) also
gives us that $\sigma_{h-l}(\lambda) = \sigma_{h-l}(\lambda')$.

We claim that $a$ must have odd height in $\lambda'$ if we are in type
$B$, and even height in types $C$ and $D$.  We prove it in type $B$ as
follows.  The first $h-l$ parts of $\lambda$ constitute a
$B$-partition, so $\sigma_{h-l}(\lambda) = \sigma_{h-l}(\lambda')$ is
odd.  Suppose $a$ had even height in $\lambda'$, and let $k$ be the
height of the next larger part of $\lambda'$.  We know that $k$ must
be even too, since $a$ must have even multiplicity.  Then
$\sigma_{h-l}(\lambda') = \sigma_k(\lambda') + (h-l-k)a$.  Since the
second term here is even, $\sigma_k(\lambda')$ must be odd.  But since
$k$ is even, the first $k$ parts of $\lambda'$ constitute a
$D$-partition, and $\sigma_k(\lambda')$ has to be even.  We have a
contradiction; therefore, $a$ has odd height in $\lambda'$.

Since $\cmp{\lambda'}{\nu'}$ is special, $a$ must have even height
with respect to $\nu'$; it cannot appear inside a basic block.  We
continue to take $k$ to be the height of the next larger part of
$\lambda'$ after $a$, but we know now that $k$ is odd.  Let $\theta' =
\chi^+_k(\lambda')$, and let $\kappa'$ be the partition consisting of
those parts of $\nu'$ that are larger than $a$.  Then
$\cmp{\theta'}{\kappa'}$ is a marked $B$-partition, and if we let $m =
r_{\lambda'}(a)$, then the expression $\cmp{\theta'}{\kappa'} \cup
\cmp{[a^m]}{\vn} \cup \cdots$ is part of a division into
blocks of $\cmp{\lambda'}{\nu'}$.  We can compute, then, that
\[
\mu'{}^* = \dS(\pmp{\theta'}{\kappa'})^* \cup [a^m] \cup \cdots.
\]
Since $\cmp{\theta'}{\kappa'}$ is of type $B$, we compute that
\[
\sigma_{h-l}(\mu'{}^*) = |\theta'|-1 + (h-l-k)a = \sigma_{h-l}(\lambda') -
1.
\]
On the other hand, in $\cmp{\lambda}{\nu}$, the part $a$ belongs to
some nonspecial basic block $\cmp{\phi}{\gamma}$, around which we can
build a division into blocks $\cmp{\theta}{\kappa} \cup
\cmp{\phi}{\gamma} \cup \cdots$.  We either have that
$\cmp{\theta}{\kappa}$ is of type $B$ and $\cmp{\phi}{\kappa}$ of type
$D$, or that $\cmp{\theta}{\kappa}$ is trivial and $\cmp{\phi}{\kappa}$
is of type $B$.  Assume we are in the former case; the latter is
handled similarly.  Let $h' = \height_\phi(a)$, and let $\phi' =
\chi^+_{h'-l}(\phi)$ and $\phi'' = \chi^-_{h'-l}(\phi)$.  
Using Proposition~\ref{prop:basic} and Lemma~\ref{lem:cupvee} to
compute $\dS(\pmp{\phi}{\gamma})$, we fine
\[
\mu^* = \dS(\pmp{\theta}{\kappa})^* \cup \phi'{}^{+*+}\!_{D-}\!^* \cup
\phi''{}^{-*}\!_B\!^* \cup \cdots.
\]
We obtain
\[
\sigma_{h-l}(\mu^*) = (|\theta'|-1) + (|\phi'|+1) =
\sigma_{h-l}(\lambda)
\]
In particular, we see that $\sigma_{h-l}(\mu^*) \not\le
\sigma_{h-l}(\mu'{}^*)$, which contradicts the assumption that $\mu^*
\le \mu'{}^*$.  Therefore, the inequalities~(\ref{eqn:sigma-compare})
hold for all $i$, and we obtain $\zeta \le \lambda'$, as desired.
\end{proof}

\begin{prop}\label{prop:ord-rev}
If $\cmp{\lambda}{\nu} \le \cmp{\lambda'}{\nu'}$, then
$\dA(\pmp{\lambda}{\nu}) \ge \dA(\pmp{\lambda'}{\nu'})$.
\end{prop}
\begin{proof}
We first prove the statement in the special case that
$\cmp{\lambda}{\nu}$ and $\cmp{\lambda'}{\nu'}$ are special.  The
$\Nocc$-inequality $\dA(\pmp{\lambda}{\nu}) \ge
\dA(\pmp{\lambda'}{\nu'})$ is equivalent to the two $\No$-inequalities
\begin{equation}\label{eqn:ord-rev}
\dS(\pmp{\lambda}{\nu}) \ge \dS(\pmp{\lambda'}{\nu'})
\qquad\text{and}\qquad
\dS(\dA(\pmp{\lambda}{\nu})) \le \dS(\dA(\pmp{\lambda'}{\nu'})).
\end{equation}
The first of these is implied by $\cmp{\lambda}{\nu} \le
\cmp{\lambda'}{\nu'}$, by definition.  For the second, since these
marked partitions are special, we know $\dS(\dA(\pmp{\lambda}{\nu})) =
\lambda$ and $\dS(\dA(\pmp{\lambda'}{\nu'})) = \lambda'$.  But the
inequality $\lambda \le \lambda'$ is again part of the definition of
$\cmp{\lambda}{\nu} \le \cmp{\lambda'}{\nu'}$.  Thus,
(\ref{eqn:ord-rev}) holds, and the proposition holds for special
marked partitions.

Now, if $\cmp{\lambda}{\nu}$ and $\cmp{\lambda'}{\nu'}$ are arbitrary
marked partitions with $\cmp{\lambda}{\nu} \le \cmp{\lambda'}{\nu'}$,
we obtain $\cmp{\lambda}{\nu} \le \dA^2(\pmp{\lambda'}{\nu'})$ by
Proposition~\ref{prop:ord-pres}.  Then, repeated application of
Lemma~\ref{lem:pre-ord-rev} implies that $e(\pmp{\lambda}{\nu}) =
\dA^2(\pmp{\lambda}{\nu}) \le \dA^2(\pmp{\lambda'}{\nu'})$. 
Both sides of this inequality are special marked partitions, so the
previous paragraph tells us that $\dA^3(\pmp{\lambda}{\nu}) \ge
\dA^3(\pmp{\lambda'}{\nu'})$.  Finally,
Proposition~\ref{prop:spec-inv}, combined with
Proposition~\ref{prop:im-special}, says that $\dA^3 = \dA$, so we get
$\dA(\pmp{\lambda}{\nu}) \ge \dA(\pmp{\lambda}{\nu})$, as desired.
\end{proof}

\subsection{Maximality of the image}

We have now established that $\dA$ satisfies each of the
axioms~(\ref{ax:ord-rev}), (\ref{ax:ord-pres}), and
(\ref{ax:dS-resp}), in Propositions~\ref{prop:ord-rev},
\ref{prop:ord-pres}, and~\ref{prop:dS-resp}, respectively.  Only
axiom~(\ref{ax:im-size}) remains, but we are not going to verify it
directly.  Instead, we employ the strategy developed in
Section~\ref{sect:formal}: we need only prove that the set of special
marked partitions corresponds to the set $\Noccsp$ defined there, and
then the theorem follows by application of
Proposition~\ref{prop:exist-uniq}. 

\begin{lem}
Let $\cmp{\lambda}{\nu}$ be a nonspecial basic block, and let $\mu =
\dS(\pmp{\lambda}{\nu})$.  There does not exist a marking partition
$\xi$ such that $\dS(\pmp{\mu}{\xi}) = \lambda$.
\end{lem}
\begin{proof}
Let us assume that $\cmp{\lambda}{\nu}$ is a marked partition of type
$B$.  Since $\cmp{\lambda}{\nu}$ is not special, we can apply the
partial specialization map to it and obtain $\cmp{\lambda'}{\nu'} =
s(\pmp{\lambda}{\nu})$, where $\lambda' > \lambda$.  A brief
consideration of the definition of $s$ reveals that
$\cmp{\lambda'}{\nu'}$ must itself be either a basic block or a
trivially marked partition; moreover, the latter can be obtained only
if $\lambda'$ has no parts of odd height except the last one.  It is
easily verified that for trivially marked partitions with this
property, $\dS$ is given by the formulas of
Proposition~\ref{prop:basic}: thus, $\dS(\pmp{\lambda'}{\nu'}) =
\lambda'{}^{-*}\!_C$ regardless of whether $\nu'$ is trivial or not.

Next, we prove that for any $\xi$,
\begin{equation}\label{eqn:max-sp}
\dS(\pmp{\mu}{\xi}) \ge \mu^{+*}\!_B \ge \lambda'.
\end{equation}
(The appropriate expressions for the middle term in types $C$ and $D$
are $\mu^{-*}\!_C$ and $\mu^{+-*}\!_D$, respectively: these
formulas are those appearing in Proposition~\ref{prop:basic}.)  The
lemma then follows, because we will have that $\dS(\pmp{\mu}{\xi}) >
\lambda$ for all $\xi$.

We shall make the assumption that the smallest part of $\lambda'$ is
$1$.  If $\nu'$ is not trivial, we are just assuming that
$\cmp{\lambda'}{\nu'}$ is an ultrabasic block; (\ref{eqn:max-sp}) then
follows for general basic blocks by Lemma~\ref{lem:square-blocks}, as
usual.  If $\nu'$ is trivial, the same reduction still works, because
the fact that only the last part of $\lambda'$ has odd height implies
that it has odd multiplicity, and is therefore odd.

For the left-hand inequality of (\ref{eqn:max-sp}), suppose
$\cmp{\mu}{\xi} = (\xi,\omega)$.  Then (note that this is a type-$C$
marked partition) we have
\[
\xi \cup \omega{}^+\!_B \le \mu^+
\]
and thence $(\xi \cup \omega^+\!_B)^*\!_B \ge \mu^{+*}\!_B$.
For the right-hand inequality, it is easily verified that
$\lambda'{}^{-*}\!_C\!^{+*} \ge \lambda'$, although we need to use the
fact that the smallest part of $\lambda'$ is $1$.  It then follows that
$\mu^{+*}\!_B = \lambda'{}^{-*}\!_C\!^{+*}\!_B \ge \lambda'_B =
\lambda'$.
\end{proof}

\begin{prop}\label{prop:max-sp-class}
Let $\cmp{\lambda}{\nu}$ be any nonspecial marked partition, and let
$\mu = \dS(\pmp{\lambda}{\nu})$.  There does not exist a marking
partition $\xi$ such that $\dS(\pmp{\mu}{\xi}) = \lambda$.
\end{prop}
\begin{proof}
We employ induction on the number of parts of $\lambda$.  The previous
lemma handles the case where $\cmp{\lambda}{\nu}$ is a basic block.
Otherwise, choose some division into blocks $\cmp{\lambda}{\nu} =
\cmp{\lambda_1}{\nu_1} \cup \cmp{\lambda_2}{\nu_2}$, and let $\mu_i =
\dS(\pmp{\lambda_i}{\nu_i})$ for $i = 1$,~$2$.  Suppose we have some
$\xi$ such that $\dS(\pmp{\mu}{\xi}) = \lambda$.  If it is possible to
write $\cmp{\mu}{\xi}$ as $\cmp{\mu_1}{\xi_1} \vee
\cmp{\mu_2}{\xi_2}$ for some $\xi_1$ and $\xi_2$, then it would follow
that $\dS(\pmp{\mu_i}{\xi_i}) = \lambda_i$ for each $i$.  But at least
one of the $\cmp{\lambda_i}{\xi_i}$ is nonspecial, so that would
contradict the inductive hypothesis.

Suppose, on the other hand, that $\cmp{\mu}{\xi}$ cannot be written as
such a join.  In this case, we again use the technique employed for
Theorem~\ref{thm:po-class} and Lemma~\ref{lem:pre-ord-rev}.  Write
$\cmp{\lambda_1}{\nu_1}$ as $\cmp{[a^l]}{\vn} \vee
\cmp{\lambda'_1}{\nu'_1}$, where the latter is a decomposition
satisfying the hypotheses of Lemma~\ref{lem:square-blocks}, chosen
such that $a$ is as large as possible.  In particular, $a$ will be at
least as large as the largest part of $\lambda_2$.  It is clear that
$\sigma_a(\lambda^*) = |\lambda_2| + al$.  A brief glance at the
formulas of Proposition~\ref{prop:blocks} shows that, moreover,
$\sigma_a(\mu) = \sigma_a(\lambda^*)$ if $\cmp{\lambda}{\nu}$ is of
type $B$ or $D$, but $\sigma_a(\mu) = \sigma_a(\lambda^*) + 1$ in type
$C$.

We now turn our attention to $\cmp{\mu}{\xi}$.  Let $\zeta =
\dS(\pmp{\mu}{\xi})$.  The assumption that $\cmp{\mu}{\xi}$ cannot be
written as an appropriate kind of join means that the $a$-th part of
$\mu$ has odd (generalized) height in $\xi$.  A calculation much like
that carried out for the proof of Theorem~\ref{thm:d-class}, whose
details we omit, shows that $\sigma_a(\zeta^*)$ is equal to
$\sigma_a(\mu)$ if $\cmp{\mu}{\xi}$ is of type $B$ ({\it i.e.}, if
$\cmp{\lambda}{\nu}$ is of type $C$), or to $\sigma_a(\mu) + 1$ if
$\cmp{\mu}{\xi}$ is of type $C$ or $D$.  Thus, in all cases, we have
$\sigma_a(\zeta^*) = \sigma_a(\lambda^*) + 1$.  In particular, this
means that $\zeta \ne \lambda$, contradicting our assumption.
\end{proof}

The preceding proposition says exactly that the set of special marked
partitions coincides with the set $\Noccsp$ defined in
Section~\ref{sect:formal}, as promised.  We therefore obtain the
following theorem.

\begin{thm}\label{thm:d-class}
The map $\dA$ is the unique extended duality map in the classical
groups. \qed 
\end{thm}

\section{Explicit calculations and the exceptional groups}
\label{sect:explicit}

The main results in the case of the exceptional groups are established
by explicit calculation.  In this section, we present explicit
calculations of the partial order and the duality map in all of the
exceptional groups, as well as in a number of classical groups of
small rank.

We name elements $(\orb,C) \in \Nocc$ in the exceptional groups by a
pair of symbols $(L_1,L_2)$, where $L_1$ is the Bala-Carter notation
for $\orb$, as found in, say, \cite{carter:finite-gps}, and $L_2$ is
the label Sommers assigns to $(\orb,C)$ in his generalized Bala-Carter
theorem \cite{sommers:b-c}.  (Of course, we are only writing down
$L_1$ for our own convenience, since $L_2$ alone determines the
orbit.)  We deviate from this notation when $C$ is the trivial
conjugacy class in $\Ab(\orb)$: in this case, the generalized
Bala-Carter label for $(\orb,C)$ is the same as the Bala-Carter label
for $\orb$, but for the sake of brevity, we write $(L_1,1)$ rather
than $(L_1,L_1)$.

A further comment about generalized Bala-Carter labels for pairs
$(\orb,C)$ is in order, because the generalized Bala-Carter theorem is
actually a classification of $\Noc$, not of $\Nocc$.  For most orbits
in the exceptional groups, we have $A(\orb) = \Ab(\orb)$, so this
distinction does not matter, but in a handful of cases, $\Ab(\orb)$
has fewer conjugacy classes than $A(\orb)$.  This occurs for two
orbits in $F_4$, two in $E_7$, and seven in $E_8$.  In all but one of
these cases, we have $A(\orb) = S_2$ and $\Ab(\orb) = 1$; however, for
the orbit $E_8(b_6)$ in type $E_8$, we have $A(\orb) = S_3$,
$\Ab(\orb) = S_2$.  In all of these cases, the only ambiguity is that
two conjugacy classes of $A(\orb)$ map to the trivial conjugacy class
of $\Ab(\orb)$.  (In the $E_8(b_6)$ example, only one conjugacy class
of $S_3$ descends to the nontrivial conjugacy class of $S_2$.)  In
each such situation, we simply ignore the nontrivial class of
$A(\orb)$ that maps to the trivial one in $\Ab(\orb)$, and we
designate the latter with a label of the form $(L_1,1)$.

\begin{thm}
Let $C, C' \subset A(\orb)$ be two conjugacy classes associated to the
same orbit.  Then, in the exceptional groups, $\dS(\orb,C) =
\dS(\orb,C')$ if and only if $C$ and $C'$ have the same image in
$\Ab(\orb)$.  As a consequence, the partial order (\ref{eqn:po}) is
well-defined in the exceptional groups.
\end{thm}
\begin{proof}
Sommers gives tables of all the values of $\dS$ on all pairs
$(\orb,C) \in \Noc$ for each exceptional group in
\cite{sommers:duality}.  We merely read through this table and verify
that the above statement is true.
\end{proof}

\begin{thm}
There exists a unique extended duality map in the case of each exceptional
group.
\end{thm}

\newcommand{\noccdiagram}[1]{\begin{tabular}{@{}c@{}}
  \scalebox{0.48}{\xymatrix@=8pt{#1}}\end{tabular}}
\newcommand{\nocclabel}[1]{\raisebox{-3pt}{\scalebox{2}{$#1$}}}
\newcommand{\emp}[2]{(#1,#2)}

\begin{table}[b]
\begin{center}
\newcommand{\co}{\vn}
\begin{tabular}{cccc}
\noccdiagram{%
*{\nocclabel{B_2}}
 &*+{\cmp{[5]}{\co}} \ar@{-}[d] \\
 &*+{\cmp{[3,1^2]}{[3,1]}} \ar@{-}[d] \\
 &*+{\cmp{[3,1^2]}{\co}} \ar@{-}[d] \\
 &*+{\cmp{[2^2,1]}{\co}} \ar@{-}[d] \\
 &*+{\cmp{[1^5]}{\co}}
\ar@{}[d] \\ &\strut\ar@{}[d] \\
*{\nocclabel{C_2}}
 &*+{\cmp{[4]}{\co}} \ar@{-}[d] \\
 &*+{\cmp{[2^2]}{[2]}} \ar@{-}[d] \\
 &*+{\cmp{[2^2]}{\co}} \ar@{-}[d] \\
 &*+{\cmp{[2,1^2]}{\co}} \ar@{-}[d] \\
 &*+{\cmp{[1^4]}{\co}}
\ar@{}[d] \\  \strut\ar@{-}[rr]&\strut\ar@{}[d]&\strut \\
*{\nocclabel{B_3}}
 &*+{\cmp{[7]}{\co}} \ar@{-}[d] \\
 &*+{\cmp{[5,1^2]}{[5,1]}} \ar@{-}[d] \\
 &*+{\cmp{[5,1^2]}{\co}} \ar@{-}[d] \\
 &*+{\cmp{[3^2,1]}{\co}} \ar@{-}[dl]\ar@{-}[dr] \\
*+{\cmp{[3,1^4]}{[3,1]}} \ar@{-}[dr] &&*+{\cmp{[3,2^2]}{\co}} \ar@{-}[dl] \\
 &*+{\cmp{[3,1^4]}{\co}} \ar@{-}[d] \\
 &*+{\cmp{[2^2,1^3]}{\co}} \ar@{-}[d] \\
 &*+{\cmp{[1^7]}{\co}}
\ar@{}[d] \\ &\strut\ar@{}[d] \\
*{\nocclabel{C_3}}
 &*+{\cmp{[6]}{\co}} \ar@{-}[d] \\
 &*+{\cmp{[4,2]}{[2]}} \ar@{-}[d] \\
 &*+{\cmp{[4,2]}{\co}} \ar@{-}[dl]\ar@{-}[dr] \\
*+{\cmp{[4,1^2]}{\co}} \ar@{-}[dr] && *+{\cmp{[3^2]}{\co}} \ar@{-}[dl] \\
 &*+{\cmp{[2^3]}{\co}} \ar@{-}[d] \\
 &*+[F]{\cmp{[2^2,1^2]}{[2]}} \ar@{-}[d] \\
 &*+{\cmp{[2^2,1^2]}{\co}} \ar@{-}[d] \\
 &*+{\cmp{[2,1^4]}{\co}} \ar@{-}[d] \\
 &*+{\cmp{[1^6]}{\co}}
}
&
\noccdiagram{
*{\nocclabel{B_4}}
 &*+{\cmp{[9]}{\co}} \ar@{-}[d] \\
 &*+{\cmp{[7,1^2]}{[7,1]}} \ar@{-}[d] \\
 &*+{\cmp{[7,1^2]}{\co}} \ar@{-}[d] \\
 &*+{\cmp{[5,3,1]}{[5,1]}} \ar@{-}[dl]\ar@{-}[d] \\
*+{\cmp{[5,1^4]}{[5,1]}} \ar@{-}[dd]
  &*+{\cmp{[5,3,1]}{\co}} \ar@{-}[d]\ar@{-}[dr] \\
 &*+{\cmp{[5,2^2]}{\co}} \ar@{-}[dl]\ar@{-}[d]
   &*+{\cmp{[4^2,1]}{\co}} \ar@{-}[dl] \\
*+{\cmp{[5,1^4]}{\co}} \ar@{-}[dr] &*+{\cmp{[3^3]}{\co}} \ar@{-}[d] \\
 &*+{\cmp{[3^2,1^3]}{\co}} \ar@{-}[d] \\
 &*+[F]{\cmp{[3,2^2,1^2]}{[3,1]}} \ar@{-}[dl]\ar@{-}[dr] \\
*+{\cmp{[3,1^6]}{[3,1]}} \ar@{-}[d]
    && *+{\cmp{[3,2^2,1^2]}{\co}} \ar@{-}[dll]\ar@{-}[d] \\
*+{\cmp{[3,1^6]}{\co}} \ar@{-}[dr]    && *+{\cmp{[2^4,1]}{\co}} \ar@{-}[dl] \\
 &*+{\cmp{[2^2,1^5]}{\co}} \ar@{-}[d] \\
 &*+{\cmp{[1^9]}{\co}}
\ar@{}[d] \\ &\strut\ar@{}[d] \\
*{\nocclabel{C_4}}
 &*+{\cmp{[8]}{\co}} \ar@{-}[d] \\
 &*+{\cmp{[6,2]}{[2]}} \ar@{-}[dl]\ar@{-}[dr] \\
*+{\cmp{[6,2]}{\co}} \ar@{-}[d]\ar@{-}[drr]
    && *+{\cmp{[4^2]}{[4]}} \ar@{-}[d] \\
*+{\cmp{[6,1^2]}{\co}} \ar@{-}[dr]          
    && *+{\cmp{[4^2]}{\co}} \ar@{-}[dl] \\
 &*+{\cmp{[4,2^2]}{\co}} \ar@{-}[dl]\ar@{-}[dd] \\
*+[F]{\cmp{[4,2,1^2]}{[2]}} \ar@{-}[d] \\
*+{\cmp{[4,2,1^2]}{\co}} \ar@{-}[dd]\ar@{-}[dr]
  &*+{\cmp{[3^2,2]}{\co}} \ar@{-}[d]\ar@{-}[dr] \\
 &*+{\cmp{[3^2,1^2]}{\co}} \ar@{-}[d] &*+{\cmp{[2^4]}{[2]}} \ar@{-}[dl] \\
*+{\cmp{[4,1^4]}{\co}} \ar@{-}[dr] &*+{\cmp{[2^4]}{\co}} \ar@{-}[d] \\
 &*+{\cmp{[2^3,1^2]}{\co}} \ar@{-}[d] \\
 &*+[F]{\cmp{[2^2,1^4]}{[2]}} \ar@{-}[d] \\
 &*+{\cmp{[2^2,1^4]}{\co}} \ar@{-}[d] \\
 &*+{\cmp{[2,1^6]}{\co}} \ar@{-}[d] \\
 &*+{\cmp{[1^8]}{\co}}
}
&
\noccdiagram{
*{\nocclabel{D_3}}
 &*+{\cmp{[5,1]}{\co}} \ar@{-}[d] \\
 &*+{\cmp{[3^2]}{\co}} \ar@{-}[d] \\
 &*+{\cmp{[3,1^3]}{\co}} \ar@{-}[d] \\
 &*+{\cmp{[2^2,1^2]}{\co}} \ar@{-}[d] \\
 &*+{\cmp{[1^6]}{\co}}
\ar@{}[d] \\ \strut\ar@{-}[rr]&\strut\ar@{}[d]&\strut \\
*{\nocclabel{D_4}}
 &*+{\cmp{[7,1]}{\co}} \ar@{-}[d] \\
 &*+{\cmp{[5,3]}{\co}} \ar@{-}[dl]\ar@{-}[d]\ar@{-}[dr] \\
*+{\cmp{[4^2]^{I}}{\co}} \ar@{-}[dr] 
  &*+{\cmp{[5,1^3]}{\co}} \ar@{-}[d]
  &*+{\cmp{[4^2]^{II}}{\co}} \ar@{-}[dl] \\
 &*+{\cmp{[3^2,1^2]}{[3,1]}} \ar@{-}[d] \\
 &*+{\cmp{[3^2,1^2]}{\co}} \ar@{-}[d] \\
 &*+{\cmp{[3,2^2,1]}{\co}} \ar@{-}[dl]\ar@{-}[d]\ar@{-}[dr] \\
*+{\cmp{[2^4]^{I}}{\co}} \ar@{-}[dr] 
  &*+{\cmp{[3,1^5]}{\co}} \ar@{-}[d] & *+{\cmp{[2^4]^{II}}{\co}} \ar@{-}[dl] \\
 &*+{\cmp{[2^2,1^4]}{\co}} \ar@{-}[d] \\
 &*+{\cmp{[1^8]}{\co}}
\ar@{}[d] \\ \strut\ar@{-}[rr]&\strut\ar@{}[d]&\strut \\
*{\nocclabel{D_5}}
 &*+{\cmp{[9,1]}{\co}} \ar@{-}[d] \\
 &*+{\cmp{[7,3]}{\co}} \ar@{-}[dl]\ar@{-}[dr] \\
*+{\cmp{[7,1^3]}{\co}} \ar@{-}[dr] &&*+{\cmp{[5^2]}{\co}} \ar@{-}[dl] \\
 &*+{\cmp{[5,3,1^2]}{[3,1]}} \ar@{-}[d] \\
 &*+{\cmp{[5,3,1^2]}{\co}} \ar@{-}[dl]\ar@{-}[dr] \\
*+{\cmp{[4^2,1^2]}{\co}} \ar@{-}[d] 
    &&*+{\cmp{[5,2^2,1]}{\co}} \ar@{-}[dll]\ar@{-}[d] \\
*+{\cmp{[3^3,1]}{\co}} \ar@{-}[d]\ar@{-}[drr] 
    &&*+{\cmp{[5,1^5]}{\co}} \ar@{-}[d] \\
*+{\cmp{[3^2,2^2]}{\co}} \ar@{-}[dr]
    &&*+{\cmp{[3^2,1^4]}{[3,1]}} \ar@{-}[dl] \\
 &*+{\cmp{[3^2,1^4]}{\co}} \ar@{-}[d] \\
 &*+{\cmp{[3,2^2,1^3]}{\co}} \ar@{-}[dl]\ar@{-}[dr] \\
*+{\cmp{[2^4,1^2]}{\co}} \ar@{-}[dr] &&*+{\cmp{[3,1^7]}{\co}} \ar@{-}[dl] \\
 &*+{\cmp{[2^2,1^6]}{\co}} \ar@{-}[d] \\
 &*+{\cmp{[1^{10}]}{\co}}
}
&
\noccdiagram{
*{\nocclabel{G_2}}
 &*+{\emp{G_2}{1}} \ar@{-}[d] \\
 &*+{\emp{G_2(a_1)}{A_2}} \ar@{-}[d] \\
 &*+{\emp{G_2(a_1)}{A_1+\tilde A_1}} \ar@{-}[d] \\
 &*+{\emp{G_2(a_1)}{1}} \ar@{-}[d]\\
 &*+{\emp{\tilde A_1}{1}} \ar@{-}[d]\\
 &*+{\emp{A_1}{1}} \ar@{-}[d]\\
 &*+{\emp{1}{1}}
\ar@{}[d] \\ &\strut\ar@{}[d] \\
\strut\ar@{-}[rr]&\strut\ar@{}[d]&\strut \\ &\strut\ar@{}[d] \\
*{\nocclabel{F_4}}
 &*+{\emp{F_4}{1}} \ar@{-}[d] \\
 &*+{\emp{F_4(a_1)}{B_4}} \ar@{-}[d] \\
 &*+{\emp{F_4(a_1)}{1}} \ar@{-}[d] \\
 &*+{\emp{F_4(a_2)}{1}} \ar@{-}[dl]\ar@{-}[dr] \\
*+{\emp{C_3}{1}} \ar@{-}[dd] &&*+{\emp{B_3}{1}} \ar@{-}[d] \\
    &&*+{\emp{F_4(a_3)}{A_3+\tilde A_1}} \ar@{-}[dll]\ar@{-}[d] \\
*+{\emp{F_4(a_3)}{A_2+\tilde A_2}} \ar@{-}[dr] 
    &&*+{\emp{F_4(a_3)}{B_4(a_1)}} \ar@{-}[dl]\ar@{-}[dd]\\
 &*+{\emp{F_4(a_3)}{A_1+C_3(a_1)}} \ar@{-}[dl]\ar@{-}[d]\\
*+{\emp{F_4(a_3)}{1}} \ar@{-}[dr]
  &*+{\emp{C_3(a_1)}{A_1+B_2}} \ar@{-}[d]
  &*+{\emp{B_2}{A_3}} \ar@{-}[dd]\\
 &*+{\emp{C_3(a_1)}{1}} \ar@{-}[dl]\ar@{-}[dr]\\
*+{\emp{\tilde A_2+A_1}{1}} \ar@{-}[dd]\ar@{-}[drr] 
    &&*+{\emp{B_2}{1}} \ar@{-}[d]\\
   &&*+{\emp{A_2+\tilde A_1}{1}} \ar@{-}[d] \\
*+{\emp{\tilde A_2}{1}} \ar@{-}[dr] &&*+{\emp{A_2}{1}} \ar@{-}[dl] \\
 &*+{\emp{A_1+\tilde A_1}{1}} \ar@{-}[d] \\
 &*+[F]{\emp{\tilde A_1}{2A_1}} \ar@{-}[d] \\
 &*+{\emp{\tilde A_1}{1}} \ar@{-}[d] \\
 &*+{\emp{A_1}{1}} \ar@{-}[d] \\
 &*+{\emp{1}{1}}
}
\end{tabular}

\end{center}
\end{table}

\begin{table}
\begin{center}
\newcommand{\co}{1}
\begin{tabular}{cc}
\noccdiagram{
*{\nocclabel{E_6}}
  &*+{\emp{E_6}{\co}} \ar@{-}[d] \\
  &*+{\emp{E_6(a_1)}{\co}} \ar@{-}[d] \\
  &*+{\emp{D_5}{\co}} \ar@{-}[d] \\
  &*+{\emp{E_6(a_3)}{A_5+A_1}} \ar@{-}[d] \\
  &*+{\emp{E_6(a_3)}{\co}} \ar@{-}[dl]\ar@{-}[d] \\
*+{\emp{A_5}{\co}} \ar@{-}[d]
  &*+{\emp{D_5(a_1)}{\co}} \ar@{-}[dl]\ar@{-}[dr] \\
*+{\emp{A_4+A_1}{\co}} \ar@{-}[d]\ar@{-}[drr]
    &&*+{\emp{D_4}{\co}} \ar@{-}[d] \\
*+{\emp{A_4}{\co}} \ar@{-}[dr] &&*+{\emp{D_4(a_1)}{3A_2}} \ar@{-}[dl] \\
  &*+{\emp{D_4(a_1)}{A_3+2A_1}} \ar@{-}[d] \\
  &*+{\emp{D_4(a_1)}{\co}} \ar@{-}[d] \\
  &*+{\emp{A_3+A_1}{\co}} \ar@{-}[dl]\ar@{-}[dr]\\
*+{\emp{A_3}{\co}} \ar@{-}[d]
    &&*+{\emp{2A_2+A_1}{\co}} \ar@{-}[dll]\ar@{-}[d] \\
*+{\emp{A_2+2A_1}{\co}} \ar@{-}[d]\ar@{-}[dr]
    &&*+{\emp{2A_2}{\co}} \ar@{-}[dl] \\
*+{\emp{A_2}{4A_1}} \ar@{-}[dr] &*+{\emp{A_2+A_1}{\co}} \ar@{-}[d] \\
  &*+{\emp{A_2}{\co}} \ar@{-}[d]\\
  &*+{\emp{3A_1}{\co}} \ar@{-}[d] \\
  &*+{\emp{2A_1}{\co}} \ar@{-}[d] \\
  &*+{\emp{A_1}{\co}} \ar@{-}[d] \\
  &*+{\emp{1}{\co}}
\ar@{}[d] \\ \strut\ar@{-}[rr]&\strut\ar@{}[d]&\strut \\
*{\nocclabel{E_7}}
  &*+{\emp{E_7}{\co}} \ar@{-}[d] \\
  &*+{\emp{E_7(a_1)}{\co}} \ar@{-}[d] \\
  &*+{\emp{E_7(a_2)}{\co}} \ar@{-}[dl]\ar@{-}[dr] \\
*+{\emp{E_6}{\co}} \ar@{-}[d]
    &&*+{\emp{E_7(a_3)}{A_1+D_6}} \ar@{-}[dll]\ar@{-}[d] \\
*+{\emp{E_6(a_1)}{A_7}} \ar@{-}[dr]
    &&*+{\emp{E_7(a_3)}{\co}} \ar@{-}[dl]\ar@{-}[d] \\
  &*+{\emp{E_6(a_1)}{\co}} \ar@{-}[d] &*+{\emp{D_6}{\co}} \ar@{-}[dl]\\
  &*+{\emp{E_7(a_4)}{\co}} \ar@{-}[dl]\ar@{-}[d]\ar@{-}[dr] \\
*+{\emp{A_6}{\co}} \ar@{-}[d]
   &*+{\emp{D_5+A_1}{\co}} \ar@{-}[dl]\ar@{-}[dr] 
   &*+{\emp{D_6(a_1)}{\co}} \ar@{-}[ddll]\ar@{-}[d] \\
*+{\emp{E_7(a_5)}{A_5+A_2}} \ar@{-}[d] &&*+{\emp{D_5}{\co}} \ar@{-}[dd] \\
*+{\emp{E_7(a_5)}{A_1+D_6(a_2)}} \ar@{-}[d]\ar@{-}[drr] \\
*+{\emp{E_7(a_5)}{\co}} \ar@{-}[d]\ar@{-}[drr]
    &&*+{\emp{E_6(a_3)}{(A_5+A_1)'}} \ar@{-}[d] \\
*+{\emp{D_6(a_2)}{\co}} \ar@{-}[d]\ar@{-}[dr]
    &&*+{\emp{E_6(a_3)}{\co}} \ar@{-}[dl]\ar@{-}[d] \\
*+{\emp{A_5+A_1}{\co}} \ar@{-}[d]\ar@{-}[drr]
   &*+{\emp{(A_5)'}{\co}} \ar@{-}[dr]
    &*+{\emp{D_5(a_1)+A_1}{\co}} \ar@{-}[dl]\ar@{-}[d]\\
*+{\emp{(A_5)''}{\co}} \ar@{-}[ddddr]
   &*+[F]{\emp{D_5(a_1)}{D_4+2A_1}} \ar@{-}[d] 
    &*+{\emp{A_4+A_2}{\co}} \ar@{-}[d] \\
  &*+{\emp{D_5(a_1)}{\co}} \ar@{-}[d]\ar@{-}[ddr]
    &*+[F]{\emp{A_4+A_1}{A_1+2A_3}} \ar@{-}[dddl]\ar@{-}[dd] \\
  &*+{\emp{D_4+A_1}{\co}} \ar@{-}[ddddl]\ar@{-}[ddddr] \\
   &&*+{\emp{A_4+A_1}{\co}} \ar@{-}[ddl]\ar@{-}[ddd] \\
  &*+{\emp{A_4}{2A_3}} \ar@{-}[d] \\
  &*+{\emp{A_4}{\co}} \ar@{-}[ddr] \\
*+{\emp{D_4}{\co}} \ar@{-}[d]
    &&*+{\emp{A_3+A_2+A_1}{\co}} \ar@{-}[dll]\ar@{-}[dl]\ar@{-}[d] \\
*+{\emp{D_4(a_1)}{3A_2}} \ar@{-}[d]
   &*+{\emp{D_4(a_1)+A_1}{A_3+3A_1}} \ar@{-}[dl]\ar@{-}[dr] 
    &*+{\emp{A_3+A_2}{\co}} \ar@{-}[d]\\
*+{\emp{D_4(a_1)}{(A_3+2A_1)'}} \ar@{-}[d]
    &&*+{\emp{D_4(a_1)+A_1}{\co}} \ar@{-}[dll]\ar@{-}[d] \\
*+{\emp{D_4(a_1)}{\co}} \ar@{-}[d]
    &&*+{\emp{A_3+2A_1}{\co}} \ar@{-}[dll]\ar@{-}[dd] \\
*+{\emp{(A_3+A_1)'}{\co}} \ar@{-}[d]\ar@{-}[ddrr] \\
*+{\emp{2A_2+A_1}{\co}} \ar@{-}[d]\ar@{-}[dr]
    &&*+{\emp{(A_3+A_1)''}{\co}} \ar@{-}[dl]\ar@{-}[d]\\
*+{\emp{A_2+3A_1}{\co}} \ar@{-}[dr]
   &*+{\emp{2A_2}{\co}} \ar@{-}[d] &*+{\emp{A_3}{\co}} \ar@{-}[dl] \\
  &*+{\emp{A_2+2A_1}{\co}} \ar@{-}[d] \\
  &*+[F]{\emp{A_2+A_1}{5A_1}} \ar@{-}[d]\ar@{-}[dr] \\
  &*+{\emp{A_2+A_1}{\co}} \ar@{-}[dl]\ar@{-}[dr]
   &*+{\emp{A_2}{(4A_1)'}} \ar@{-}[d] \\
*+{\emp{4A_1}{\co}} \ar@{-}[d]\ar@{-}[drr] &&*+{\emp{A_2}{\co}} \ar@{-}[d] \\
*+{\emp{(3A_1)''}{\co}} \ar@{-}[dr] &&*+{\emp{(3A_1)'}{\co}} \ar@{-}[dl] \\
  &*+{\emp{2A_1}{\co}} \ar@{-}[d] \\
  &*+{\emp{A_1}{\co}} \ar@{-}[d] \\
  &*+{\emp{1}{\co}}
}
&
\noccdiagram{
*{\nocclabel{E_8}}
  &&*+{\emp{E_8}{\co}} \ar@{-}[d] \\
  &&*+{\emp{E_8(a_1)}{\co}} \ar@{-}[d] \\
  &&*+{\emp{E_8(a_2)}{\co}} \ar@{-}[d] \\
  &&*+{\emp{E_8(a_3)}{E_7+A_1}} \ar@{-}[dl]\ar@{-}[d] \\
 &*+{\emp{E_8(a_4)}{D_8}} \ar@{-}[dr]
  &*+{\emp{E_8(a_3)}{\co}} \ar@{-}[dl]\ar@{-}[d] \\
 &*+{\emp{E_7}{\co}} \ar@{-}[dr] &*+{\emp{E_8(a_4)}{\co}} \ar@{-}[d] \\
  &&*+{\emp{E_8(b_4)}{\co}} \ar@{-}[dl]\ar@{-}[dr] \\
 &*+{\emp{E_7(a_1)}{\co}} \ar@{-}[dddd]
   &&*+{\emp{E_8(a_5)}{D_8(a_1)}} \ar@{-}[dr] \\
    &&&&*+{\emp{E_8(a_5)}{\co}} \ar@{-}[dl]\ar@{-}[d] \\
   &&&*+{\emp{E_8(b_5)}{E_6+A_2}} \ar@{-}[dl]\ar@{-}[d]
    &*+{\emp{D_7}{\co}} \ar@{-}[dl] \\
  &&*+{\emp{E_8(b_5)}{E_7(a_2)+A_1}} \ar@{-}[dl]\ar@{-}[d]
   &*+{\emp{E_8(a_6)}{A_8}} \ar@{-}[dl] \\
 &*+{\emp{E_8(b_5)}{\co}} \ar@{-}[dl]\ar@{-}[dr]
  &*+{\emp{E_8(a_6)}{D_8(a_2)}} \ar@{-}[d] \\
*+{\emp{E_7(a_2)}{\co}} \ar@{-}[d]\ar@{-}[dr]
  &&*+{\emp{E_8(a_6)}{\co}} \ar@{-}[dl]\\
*+{\emp{E_6+A_1}{\co}} \ar@{-}[d]
 &*+{\emp{D_7(a_1)}{\co}} \ar@{-}[d]\ar@{-}[dr]\\
*+{\emp{E_6}{\co}} \ar@{-}[dddr]
  &*+[F]{\emp{E_7(a_3)}{A_1+D_6}} \ar@{-}[d]\ar@{-}[dr] 
  &*+{\emp{E_8(b_6)}{D_8(a_3)}} \ar@{-}[d]\ar@{-}[dr] \\
 &*+{\emp{E_7(a_3)}{\co}} \ar@{-}[ddr]\ar@{-}[dr] 
  &*+[F]{\emp{E_6(a_1)+A_1}{A_7+A_1}} \ar@{-}[ddl]\ar@{-}[d]
   &*+{\emp{E_8(b_6)}{\co}} \ar@{-}[dl]\ar@{-}[dr] \\
  &&*+{\emp{E_6(a_1)+A_1}{\co}} \ar@{-}[ddl]\ar@{-}[ddrr]
    &&*+{\emp{A_7}{\co}} \ar@{-}[d]\\
 &*+{\emp{E_6(a_1)}{(A_7)''}} \ar@{-}[d]
  &*+{\emp{D_6}{\co}} \ar@{-}[ddr]
    &&*+{\emp{D_7(a_2)}{D_5+A_3}} \ar@{-}[d] \\
 &*+{\emp{E_6(a_1)}{\co}} \ar@{-}[dd] &&&*+{\emp{D_7(a_2)}{\co}} \ar@{-}[dl] \\
   &&&*+{\emp{D_5+A_2}{\co}} \ar@{-}[dll]\ar@{-}[ddll]\ar@{-}[d] \\
 &*+{\emp{E_7(a_4)}{\co}} \ar@{-}[ddl]\ar@{-}[drr]
   &&*+{\emp{A_6+A_1}{\co}} \ar@{-}[dl]\ar@{-}[d] \\
 &*+{\emp{D_6(a_1)}{D_5+2A_1}} \ar@{-}[dl]\ar@{-}[d]
  &*+{\emp{E_8(a_7)}{2A_4}} \ar@{-}[dl]\ar@{-}[dr]
   &*+{\emp{A_6}{\co}} \ar@{-}[d]\\
*+{\emp{D_6(a_1)}{\co}} \ar@{-}[d]\ar@{-}[dr]
 &*+{\emp{E_8(a_7)}{A_5+A_2+A_1}} \ar@{-}[d]\ar@{-}[drr] 
   &&*+{\emp{E_8(a_7)}{D_5(a_1)+A_3}} \ar@{-}[dll]\ar@{-}[d] \\
*+{\emp{D_5+A_1}{\co}} \ar@{-}[dd]\ar@{-}[dr]
 &*+{\emp{E_8(a_7)}{E_6(a_3)+A_2}} \ar@{-}[d]\ar@{-}[dr] 
                  &&*+{\emp{E_8(a_7)}{D_8(a_5)}} \ar@{-}[dl]\ar@{-}[dd] \\
 &*+{\emp{E_7(a_5)}{A_5+A_2}} \ar@{-}[d]
                  &*+{\emp{E_8(a_7)}{E_7(a_5)+A_1}} \ar@{-}[dl]\ar@{-}[d] \\
*+{\emp{D_5}{\co}} \ar@{-}[dd]
 &*+{\emp{E_7(a_5)}{A_1+D_6(a_2)}} \ar@{-}[d]\ar@{-}[dr] 
                  &*+{\emp{E_8(a_7)}{\co}} \ar@{-}[d] 
                  &*+{\emp{D_6(a_2)}{D_4+A_3}} \ar@{-}[dd]\\
 &*+{\emp{E_6(a_3)+A_1}{A_5+2A_1}} \ar@{-}[dl]\ar@{-}[d]
                  &*+{\emp{E_7(a_5)}{\co}} \ar@{-}[dl]\ar@{-}[dr]\\
*+{\emp{E_6(a_3)}{(A_5+A_1)''}} \ar@{-}[d] 
 &*+{\emp{E_6(a_3)+A_1}{\co}} \ar@{-}[dl]\ar@{-}[d]\ar@{-}[drr]
                  &&*+{\emp{D_6(a_2)}{\co}} \ar@{-}[dll]\ar@{-}[d]\\
*+{\emp{E_6(a_3)}{\co}} \ar@{-}[ddr]\ar@{-}[dr]
 &*+{\emp{A_5+A_1}{\co}} \ar@{-}[d]\ar@{-}[dr] 
                  &&*+{\emp{D_5(a_1)+A_2}{\co}} \ar@{-}[dl]\ar@{-}[d] \\
 &*+{\emp{A_5}{\co}} \ar@{-}[ddrr]&*+{\emp{A_4+A_3}{\co}} \ar@{-}[dr]
                  &*+{\emp{D_4+A_2}{\co}} \ar@{-}[dll]\ar@{-}[d] \\
 &*+{\emp{D_5(a_1)+A_1}{\co}} \ar@{-}[d]\ar@{-}[drr]
  &&*+{\emp{A_4+A_2+A_1}{\co}} \ar@{-}[d]\\
 &*+[F]{\emp{D_5(a_1)}{D_4+2A_1}} \ar@{-}[dd] 
  &&*+{\emp{A_4+A_2}{\co}} \ar@{-}[dddl]\ar@{-}[dr] \\
   &&&&*+[F]{\emp{A_4+2A_1}{D_4(a_1)+A_3}} \ar@{-}[dl]\ar@{-}[d] \\
 &*+{\emp{D_5(a_1)}{\co}} \ar@{-}[d]\ar@{-}[ddr] 
                  &&*+[F]{\emp{A_4+A_1}{A_1+2A_3}} \ar@{-}[ddl]
                  &*+{\emp{A_4+2A_1}{\co}} \ar@{-}[ddll]\ar@{-}[d] \\
 &*+{\emp{D_4+A_1}{\co}} \ar@{-}[dddl]\ar@{-}[dddr] 
  &*+{\emp{A_4}{(2A_3)''}} \ar@{-}[ddl]
   &&*+{\emp{2A_3}{\co}} \ar@{-}[d] \\
  &&*+{\emp{A_4+A_1}{\co}} \ar@{-}[dl]\ar@{-}[dr]
                  &&*+{\emp{D_4(a_1)+A_2}{A_3+A_2+2A_1}} \ar@{-}[dl] \\
 &*+{\emp{A_4}{\co}} \ar@{-}[dd] 
                  &&*+{\emp{D_4(a_1)+A_2}{\co}} \ar@{-}[dl] \\
*+{\emp{D_4}{\co}} \ar@{-}[d]
   &&*+{\emp{A_3+A_2+A_1}{\co}} \ar@{-}[dl]\ar@{-}[d] \\
*+{\emp{D_4(a_1)}{3A_2}} \ar@{-}[dd] &*+{\emp{A_3+A_2}{\co}} \ar@{-}[d] 
  &*+[F]{\emp{D_4(a_1)+A_1}{3A_2+A_1}} \ar@{-}[dl]\\
 &*+[F]{\emp{D_4(a_1)+A_1}{A_3+3A_1}} \ar@{-}[dl]\ar@{-}[dr] \\
*+{\emp{D_4(a_1)}{(A_3+2A_1)''}} \ar@{-}[dr]
  &&*+{\emp{D_4(a_1)+A_1}{\co}} \ar@{-}[dl]\ar@{-}[d] \\
 &*+{\emp{D_4(a_1)}{\co}} \ar@{-}[dr]\ar@{-}[dddd]
  &*+{\emp{A_3+2A_1}{\co}} \ar@{-}[d]\ar@{-}[dr] \\
  &&*+{\emp{A_3+A_1}{\co}} \ar@{-}[dr]
   &*+{\emp{2A_2+2A_1}{\co}} \ar@{-}[d]\ar@{-}[dr] \\
   &&&*+{\emp{2A_2+A_1}{\co}} \ar@{-}[dr]
    &*+{\emp{2A_2}{A_2+4A_1}} \ar@{-}[d] \\
  &&&&*+{\emp{2A_2}{\co}} \ar@{-}[dl] \\
 &*+{\emp{A_3}{\co}} \ar@{-}[dr] &&*+{\emp{A_2+3A_1}{\co}} \ar@{-}[dl] \\
  &&*+{\emp{A_2+2A_1}{\co}} \ar@{-}[d] \\
  &&*+[F]{\emp{A_2+A_1}{5A_1}} \ar@{-}[dl]\ar@{-}[d] \\
 &*+{\emp{A_2}{(4A_1)''}} \ar@{-}[dr]
  &*+{\emp{A_2+A_1}{\co}} \ar@{-}[dl]\ar@{-}[d] \\
 &*+{\emp{4A_1}{\co}} \ar@{-}[dr] &*+{\emp{A_2}{\co}} \ar@{-}[d] \\
  &&*+{\emp{3A_1}{\co}} \ar@{-}[d] \\
  &&*+{\emp{2A_1}{\co}} \ar@{-}[d] \\
  &&*+{\emp{A_1}{\co}} \ar@{-}[d] \\
  &&*+{\emp{1}{\co}}
}
\end{tabular}

\end{center}
\end{table}

\begin{proof}
Once we have drawn out the partial-order diagram of $\Nocc$ for the
exceptional groups, we produce the extended duality map by working
backwards from the results of Section~\ref{sect:formal}.  Recall,
from that section, the definition of $\Noccsp$: this ought to be the
special set for the extended duality map.  We verify by inspection in
each type that for each pair $(\orb,C) \notin \Noccsp$, there is a
unique smallest element of $\Noccsp$ that is larger than it.  Next, we
define the map $\db$ by referring to Proposition~\ref{prop:d3-d} and
the proof of Proposition~\ref{prop:spec-uniq}: the latter tells us how
to compute $\db$ on special pairs, while the former does the same for
nonspecial pairs.  Finally, we tediously verify that the map thus
produced does, in fact, satisfy the first three axioms for an extended
duality map.  The theorem follows by application of
Proposition~\ref{prop:exist-uniq}.
\end{proof}

Below, we have drawn out the full Hasse diagram of the partial-order
structure on $\Nocc$ in types $B$ and $C$ up to rank $4$, in type $D$
up to rank $5$, and in all the exceptional groups.  In these diagrams,
most pairs $(\orb,C)$ are special.  Ones that are not special are
indicated by a solid box $\xymatrix{*+[F]{\ \ }}$.  The number of
elements in $\Nocc$ for $G_2$ (resp.~$F_4$, $E_6$, $E_7$, $E_8$) is
$7$ (resp.~$24$, $25$, $58$, $106$), and the number of special pairs
is $7$ (resp.~$23$, $25$, $55$, $98$).

In type $D$ and the exceptional groups, the duality map $\dA$ itself
can be visualized as follows: if the nonspecial pairs are deleted from
the diagram, the remaining partial-order diagram has a horizontal axis
of symmetry.  The duality map on special pairs is given by reflection
across this axis; then, Proposition~\ref{prop:d3-d} tells us how to
compute $\dA$ on nonspecial pairs.  For types $B$ and $C$, we have
drawn the Hasse diagram of $\Nocc(B_n)$ directly above that of
$\Nocc(C_n)$.  This combined diagram has a horizontal axis of symmetry
if nonspecial pairs are deleted, and $\dA$ is given by reflecting
across that.

The observant reader may remark upon an apparent discrepancy
between our diagram for $F_4$ and those given in other sources (such as
\cite{spaltenstein:classes} or \cite{carter:finite-gps}) for the
classical duality map.  Those sources show that the dual orbit to
$B_3$ is $\tilde{A}_2$, while the dual of $C_3$ is $A_2$.  In our
diagram, it looks as though $B_3$ and $C_3$ have been exchanged.  In
fact, this discrepancy arises because those sources are illustrating
the map $\dLS$, whereas $\db$ satisfies a compatibility condition with
$\dBV$.  Remarkably, $\dLS$ and $\dBV$ do not coincide for $F_4$, even
though it is isomorphic to its Langlands dual.  The reason is that
passing to the Langlands dual exchanges the long and short roots of
the root system, so in corresponding representations of the Weyl
groups, the action of the simple reflections corresponding to long and
short roots must be interchanged.  Alvis \cite{alvis:tables} describes
this corespondence explicitly.

\section{Further comments}
\label{sect:comments}

This final section is devoted to exploring how the partial order and
the extended duality map for $\Nocc$ can be employed to enhance
understanding and further the study of a handful of topics.  I am
especially indebted to E.~Sommers and A.-M.~Aubert for discussions
about these matters.  Some of the ideas and assertions in this section
are the product of joint work with A.-M.~Aubert, and will be properly
developed and proved in a forthcoming joint paper
\cite{achar-aubert:fc-pre}.

In Section~\ref{subsect:coxeter}, we show how our new tools can be
used to give a uniform approach to existing disparate descriptions of
the structure of the groups $\Ab(\orb)$.  In
Section~\ref{subsect:equiv-k-thy}, we revisit a conjecture made in
\cite{achar-sommers:local-sys} regarding the equivariant $K$-theory of
the nilpotent cone.  We show how to restate it in the language of the
partial order on $\Nocc$, and we then investigate a refinement of the
conjecture suggested by the rephrasing.  Finally, in
Section~\ref{subsect:springer}, we consider what the partial order
might have to say about representations of Weyl groups, via the
Springer correspondence.

One issue that we will not address, however, is that of giving an
``intrinsic'' construction of the duality map.  That is, it would be
nice to have some representation-theoretic construction explaining why
$\db(\orb,C)$ should be associated to $(\orb,C)$, rather than merely
an opaque existential statement regarding the set $\Noccsp$.  Such a
construction is likely to elucidate many aspects of the duality map.
For instance, what makes a pair $(\orb,C)$ special or nonspecial?  Is
there a way to realize $\Nocc$ geometrically, identifying its elements
with certain locally closed subvarieties of some variety in such a way
that its partial order just becomes the usual closure order?  Finally,
is there an analogue of the compatibility that $\dLS$ and $\dBV$ enjoy
with induction of nilpotent orbits?  An answer to this last question
would, of course, require a theory of induction for $\Nocc$.

\subsection{The $\boldsymbol{\Ab(\orb)}$ groups as Coxeter groups}
\label{subsect:coxeter}

In the exceptional groups, $\Ab(\orb)$ is always just a symmetric
group, and therefore a Coxeter group.  Indeed, it has a unique
Coxeter presentation up to conjugacy.  This structure was employed by
Lusztig \cite{lusztig:notes} to obtain a correspondence between
conjugacy classes of $\Ab(\orb)$ on the one hand, and parabolic
subgroups on the other.
In the classical groups, however, $\Ab(\orb)$ is a product of many
copies of $\Z/2\Z$: regarding this as a $\Z/2\Z$-vector space, any
basis is a set of simple reflections for a Coxeter presentation.
Moreover, no two such presentations are even conjugate.  Carrying out
an analogue of the constructions in \cite{lusztig:notes} requires
choosing a particular Coxeter presentation.  This is done for the
classical groups in \cite{achar-sommers:local-sys} by choosing the
simple reflections to be elements of those nontrivial conjugacy
classes $C$ for which $\dS(\orb,C)$ has maximal dimension.  That
turned out to be the correct choice for a certain conjecture regarding
local systems, which will be discussed in the next section.

The partial order on $\Nocc$ can be used to give a uniform description
of the canonical Coxeter structure of $\Ab(\orb)$ in all types.  If we
restrict our attention to a single orbit $\orb$, then the conjugacy
classes of $\Ab(\orb)$ inherit a partial order from $\Nocc$.  The
trivial conjugacy class is the smallest element in this partial order,
according to Proposition~\ref{prop:triv-min}.  Let us call a class $C$
\emph{superminimal} if it lies just above the trivial class: that is,
if $C > C'$ implies that $C'$ is the trivial class.  The choice
of simple reflections in \cite{achar-sommers:local-sys} consists
precisely of elements of superminimal conjugacy classes.  The
following result will be proved in \cite{achar-aubert:fc-pre}; it has
also been independently obtained by Sommers \cite{sommers:pc}.

\begin{thm}\label{thm:coxeter}
There is a set of involutions $S \subset \Ab(\orb)$, unique up to
conjugacy, such that:
\begin{enumgen}{1}
\item every element of $S$ is a member of a superminimal conjugacy class,
\item every superminimal conjugacy class has at least one
representative in $S$, and
\item $S$ constitutes a set of simple reflections for a presentation
of $\Ab(\orb)$ as a Coxeter group.
\end{enumgen}
\end{thm}

When studying the representations of Coxeter groups, we have available
to us the Macdonald-Lusztig-Spaltenstein operation of ``truncated
induction'' or ``$j$-induction.''  This operation is defined for a
certain class of irreducible representations, which includes all
special representations.  The truncated induction of an irreducible
representation (when it is defined) is the unique irreducible
component of the induced representation that occurs in as small a
symmetric power of the reflection representation as possible.  It
turns out that every representation of $\Ab(\orb)$ arises as the
truncated induction of the sign representation of some parabolic
subgroup, which is uniquely determined up to conjugacy.  Parabolic
subgroups are, in turn, determined by subsets of the set of simple
reflections.

Let $S$ be a set of simple reflections as found by
Theorem~\ref{thm:coxeter}, and let $P \subset S$ be a subset.  We thus
associate a certain representation $\rho_P$ of $\Ab(\orb)$ to $P$, and
we take $C_P$ to be the conjugacy class containing the product of all
the elements of $P$.  (In \cite{achar-aubert:fc-pre}, for technical
reasons, the formula for $\rho_P$ is not simply the truncated
induction of the sign representation, but rather that tensored with
the sign representation of $\Ab(\orb)$.)  The following proposition,
relating conjugacy classes and representations of $\Ab(\orb)$,
collects and rephrases facts that are implicit in the work of Lusztig
\cite{lusztig:notes} for the exceptional groups, and in
\cite{achar-sommers:local-sys} for the classical groups.

\begin{prop}\label{prop:irr-con-bij}
$C_P$ is well-defined, {\it i.e.}, independent of the order in which
the elements of $P$ are written.  Moreover, every conjugacy class of
$\Ab(\orb)$ occurs as some $C_P$, where $P$ is uniquely determined up
to conjugacy.  Therefore, the map
\[
\rho_P \leftrightsquigarrow  C_P
\]
is a natural bijection between irreducible representations and
conjugacy classes of $\Ab(\orb)$.  In addition, we have that $(\orb,C_P)
\le (\orb,C_Q)$ if and only if $P$ is conjugate to a subset of $Q$.
\end{prop}

Lusztig uses the correspondence between conjugacy classes and
parabolic subgroups to study a certain map assigning to each nilpotent
orbit an element of $\Nocc$.  We now recall the construction of that
map, and we consider what can be said about it with the aid of the
partial order.  Recall that a \emph{special piece} is the union of a
special orbit and all orbits in its closure that are not contained in
the closure of any other special orbit.  Let $\orb$ be a special
orbit, and define $\mathcal{M}(\Ab(\orb))$ to be the set of
$\Ab(\orb)$-conjugacy classes of pairs $(x,\rho)$, where $x \in
\Ab(\orb)$ and $\rho$ is an irreducible representation of the
centralizer of $x$ in $\Ab(\orb)$.  There is a natural imbedding of
the set of representations in the two-sided cell of the Weyl group
that corresponds to $\orb$ into the set $\mathcal{M}(\Ab(\orb))$.

Lusztig's map associates each orbit in the special piece containing
$\orb$ to some pair $(\orb,C)$, by examining the image of the Springer
representation of the given orbit under the above imbedding (see
\cite{lusztig:characters} and \cite{lusztig:notes}).  In particular,
$\orb$ itself is sent to $(\orb,1)$.  One proposition for exceptional
groups that appears in \cite{lusztig:notes} is equivalent to the
following tidy statement in terms of the partial order on $\Nocc$.

\begin{prop}\label{prop:lusztig-order}
Let $\orb_1$, $\orb_2$ be two nilpotent orbits in the same special
piece, assigned to $(\orb,C_1)$, $(\orb,C_2)$, respectively, by
Lusztig's map.  Then $\orb_1 \subset \overline{\orb_2}$ if and only if
$(\orb,C_1) \ge (\orb,C_2)$.
\end{prop}

This is proved in \cite{lusztig:notes} by case-by-case computation,
but an appropriate application of the extended duality map renders
this proposition obvious, as follows.  Sommers' canonical inverse (see
Section~\ref{sect:formal}) is conjectured to coincide with Lusztig's
map, under an appropriate identification $\Ab(\orb) \simeq
\Ab(\dBV(\orb))$.  Sommers has verified this conjecture for the
exceptional groups.  The identification of conjugacy classes in
$\Ab(\orb)$ with those in $\Ab(\dBV(\orb))$ is order-preserving,
because it respects their Coxeter structures.  The canonical inverse
map itself is order-reversing (Proposition~\ref{prop:dS-inv}), so
Lusztig's map is order-reversing as well.

Proposition~\ref{prop:lusztig-order} is only stated by Lusztig for the
exceptional groups because he did not have the correspondence between
conjugacy classes and parabolic subgroups of $\Ab(\orb)$ available for
classical groups, but the proposition should be true in general.

\subsection{Equivariant $\boldsymbol{K}$-theory of the nilpotent cone}
\label{subsect:equiv-k-thy}

Let $\Nor$ be the set of pairs $\{(\orb,\rho)\}$, where $\orb$ is a
nilpotent orbit in $\Lie{g}$, and $\rho$ is an irreducible
representation of the isotropy group of $\orb$ in $G$.  Let $\Noro$
(resp.~$\Norro$) be similarly defined, except that we take $\rho$ to
be a representation of $A(\orb)$ (resp.~$\Ab(\orb)$) instead.  There
are obvious inclusions $\Norro \hookrightarrow \Noro \hookrightarrow
\Nor$, given by pulling back representations.
Proposition~\ref{prop:irr-con-bij} yields a natural bijection between
$\Nocc$ and $\Norro$.  In this section, we consider those sets to be
identified; we freely use the partial order, as well as terms like
``special,'' in reference to elements of $\Norro$.

Lusztig and Vogan have independently conjectured the existence of a
bijection between $\Nor$ and the set $\DomWeights$ of dominant weights
of $G$, that should arise by studying the equivariant $K$-theory of
the nilpotent cone.  This idea has been investigated by Bezrukavnikov
\cite{bezrukavnikov:tensor}, \cite{bezrukavnikov:quasi-exc}, Ostrik
\cite{ostrik:equiv-k-theory}, and the author \cite{achar:thesis}.  In
\cite{achar:thesis}, the bijection is established for $GL(n)$ by an
explicit combinatorial algorithm.  In \cite{bezrukavnikov:quasi-exc},
the bijection is proved in general by a study of perverse equivariant
coherent sheaves on the nilpotent cone.

Now, nilpotent orbits in $\Langdual\Lie{g}$ are labelled by their
weighted Dynkin diagrams, which may be regarded as weights for $G$.
(The weighted Dynkin diagram of the orbit is the semisimple element of
the Jacobson-Morozov $\sln(2)$-triple for the orbit.)  It has been
observed that, given an orbit $\orb \in \Langdual\No$, this bijection
often sends its weighted Dynkin diagram to some pair $(\orb',\rho)$,
where $\orb' = \dBV(\orb)$ and, moreover, $\rho$ is a representation
that descends to the group $A(\orb)$.  This is mentioned in
\cite{chm-ostrik:dist-inv}; a more thorough discussion can be found in
Section~3 of \cite{achar-sommers:local-sys}.  In that paper, a
specific conjecture about $\rho$ was made.

We shall now review the conjecture of \cite{achar-sommers:local-sys},
and examine how to reformulate it using the new language of duality
for $\Nocc$ and $\Norro$.  Starting with an orbit $\orb$ in a
classical group, let $(\orb',C)$ be Sommers' canonical inverse for it.
To $C$ one associates a subgroup $H_C \subset \Ab(\orb')$, by first
giving a specific presentation of $\Ab(\orb')$ as a Coxeter group,
then expressing an element of $C$ as a product of certain simple
reflections, and finally taking $H_C$ to be the subgroup generated by
those simple reflections.  Conjecture~3.1 of
\cite{achar-sommers:local-sys} says that the Dynkin diagram of $\orb$
is associated by Lusztig's bijection to a pair $(\orb',\rho)$, where
$\rho$ is a representation occurring in $\Ind_{\scriptscriptstyle
H_C}^{\scriptscriptstyle \Ab(\orb)} 1$.

As remarked in Section~\ref{subsect:coxeter}, the presentation of
$\Ab(\orb')$ chosen in \cite{achar-sommers:local-sys} is the same as
that produced by Theorem~\ref{thm:coxeter}.  Let us identify $C$ as
some $C_P$, following Proposition~\ref{prop:irr-con-bij}.  Now, $\rho$
in turn is equal to $\rho_Q$ for some subset $Q \subset S$.  This
representation occurs in $\Ind_{\scriptscriptstyle
H_C}^{\scriptscriptstyle \Ab(\orb)} 1$ if and only if the trivial
representation occurs in the restriction of $\rho_Q$ to $H_C$, by
Frobenius reciprocity.  Moreover, the definition of $\rho_Q$ turns out
to have the consequence that the trivial representation occurs in its
restriction to the subgroup $H_C$ if and only if $P \subset Q$.  Using
Proposition~\ref{prop:irr-con-bij} again, we obtain the following
equivalent statement.

\begin{conj}[\cite{achar-sommers:local-sys}, Conjecture~3.1]
\label{conj:dynkin}
The Dynkin diagram of $\orb$ is assigned to a pair $(\orb',\rho)$ such
that $\orb' = \dBV(\orb)$ and $(\orb',\rho) \ge \db(\orb,1)$.
\end{conj}

\noindent
Sommers (\cite{achar-sommers:local-sys}, Remark~3.2) has also shown by
example that the above inequality can, indeed, fail to be an equality.
Namely, if one takes $\orb$ to be the subregular orbit in $B_n$, then
$\orb'$ is an orbit in $C_n$ with $\Ab(\orb') \simeq
\Z/2\Z$.  Let $\epsilon$ denote the nontrivial representation of
$\Ab(\orb')$.  It turns out that $\db(\orb,1) = (\orb',\epsilon)$, but
the Dynkin diagram of $\orb$ is assigned to $(\orb',1)$ when $n$ is
odd and $(\orb',\epsilon)$ when $n$ is even.  Notably, $(\orb',1)$ is
not special for $n \ge 3$.  Computed examples suggest that this may be
a necessary condition for the inequality above to fail to be an
equality.  The above conjecture can therefore be refined as follows.

\begin{genthm}{Conjecture \ref{conj:dynkin}$\boldsymbol{'}$}
The Dynkin diagram of $\orb$ is assigned to a pair $(\orb',\rho)$ such
that $\orb' = \dBV(\orb)$ and $(\orb',\rho) \ge \db(\orb,1)$, with
equality if all pairs $(\orb',\rho') > \db(\orb,1)$ are special.
\end{genthm}

\subsection{The Springer correspondence}
\label{subsect:springer}

Once we have a partial-order structure for certain local systems on
nilpotent orbits, an intriguing avenue of inquiry is the relationship
of this structure to representations of the Weyl group, via the
Springer correspondence.  Of course, the Springer correspondence
relates Weyl group representations to elements of $\Noro$, but we only
have a partial order structure on $\Norro$.  Nevertheless, we shall
put aside this stumbling block for the moment.

The statements below will actually be made in the context of the
generalized Springer correspondence, which we now review.  Let $W$ be
the Weyl group, let $u \in \orb$ be a nilpotent element, and let
$\mathcal{B}_u$ be the variety of Borel subalgebras containing $u$.
The original Springer correspondence was obtained by defining an
action of $W$ on the top-dimensional cohomology of $\mathcal{B}_u$.
It turns out that the map $\Irr(W) \to \Noro$ is injective but not, in
general, surjective.  Lusztig \cite{lusztig:ic-cplx} extended the
correspondence to account for the missing elements of $\Noro$.  In
this generalized version, the missing elements correspond to
irreducible representations of certain groups $W_L^G = N_G(L)/L$,
where $L$ is a ``cuspidal'' Levi subgroup and $N_G(L)$ is its
normalizer.  Let
\[
\nu : \coprod \Irr(W_L^G) \xrightarrow{\sim} \Noro
\]
be the bijection obtained in this way.

We now introduce a certain class of subgroups for Weyl groups that
will be required for the subsequent discussion.  Let $S$ be a set of
simple reflections generating $W$.  Furthermore, let $s_0$ be the
reflection corresponding to the highest root in the root system for
$W$, and let $S_0 = S \cup \{s_0\}$.  Now, a subgroup generated by a
subset of $S$ is called a parabolic subgroup.  Let us call a subgroup
generated by a proper subset of $S_0$ a \emph{pseudoparabolic
subgroup}.  

(This is by analogy with Sommers' term \emph{pseudo-Levi} for a
connected reductive subgroup of $G$ corresponding to the root system
generated by a given proper subset of $S_0$.  This seems to be a
synonym for \emph{endoscopic subgroup}, although that term is
unappealing when one is not doing any endoscopy theory.  The idea for
this class of subgroups has, at any rate, been in use for much longer
than Sommers' terminology: Spaltenstein \cite{spaltenstein:classes},
for instance, employs them without giving them any name whatsoever.)

The following desideratum for the relationship between representations
of the $W_L^G$ and the partial order was originally suggested by
Aubert.

\begin{des}\label{des:springer}
Let $L$ be a cuspidal Levi subgroup of $G$, let $W = W^G_L$, and let
$W'$ be any pseudoparabolic subgroup of $W$.  For any irreducible
representation $\rho$ of $W'$, there is a unique irreducible
representation $\pi$ occurring in $\Ind_{W'}^W \rho$ such that
$\nu(\pi) \ge \nu(\pi')$ for all irreducible representations $\pi'$
occurring in $\Ind_{W'}^W \rho$.  Moreover, $\pi$ occurs with
multiplicity $1$, and it coincides with the truncated induction
$j_{W'}^W \rho$.
\end{des}

One application of this statement will be the strengthening
of known results on the unipotent supports of character sheaves: this
is the principal topic of investigation in
\cite{achar-aubert:fc-pre}.  In the absence of such a statement,
previous treatments of this topic have often relied on assumptions
about dimension.  For example, in Section~4 of
\cite{geck:char-sheaves}, Geck defines a certain class of special
representations of pseudoparabolic subgroups, for which it is assumed
that $j(\rho)$ is attached to an orbit of larger dimension than any
other term of $\Ind \rho$, and then establishes a number of results
under the assumption that one is only dealing with special
representations from this class.  Similarly, Lusztig, in Theorem~10.7
of \cite{lusztig:unip-supp}, proves a statement asserting the
existence of a unique unipotent class of maximal dimension having
certain properties.  Both of these developments rely on the Springer
correspondence, so it seems likely that revisiting them with the help
of the above desideratum would lead to a considerable sharpening of
the results obtained.  In particular, Geck gives an example
(\cite{geck:char-sheaves}, Example~6.4) showing what can go wrong with
representations not belonging to his class.
Desideratum~\ref{des:springer} ought to allow a rephrasing of his
results that would accomodate such examples.

All this discussion is, of course, moot if we do not actually have a
partial order on $\Noro$: this is the stumbling block that we put
aside earlier.  In \cite{achar-aubert:fc-pre}, it will be shown how to
construct a map $\Noro \twoheadrightarrow \Norro$ that lets one pull
back the partial order.  Of course, in some respects, the partial
order on $\Noro$ cannot be as nice as that on $\Norro$: for instance,
the sign and trivial representations of a given $A(\orb)$ might fail
to be comparable in $\Noro$, whereas for $\Ab(\orb)$, the sign
representation is always smaller than the trivial one (see
Proposition~\ref{prop:irr-con-bij}).  This partial order will,
however, turn out to satisfy the above desideratum, with corresponding
implications for the study of character sheaves.

\bibliographystyle{pnaplain}
\bibliography{pramod}
\end{document}